\RequirePackage{fix-cm}

\documentclass[smallextended,numbook,runningheads]{svjour3}       

\smartqed  
\usepackage{graphicx}

\usepackage{epsfig,amsmath,amssymb,latexsym,amsfonts,curves,eepic,epsfig,epic,graphics,graphicx}
\usepackage{hyperref}
\usepackage{trees}

\begin{document}

\title{A sixth order averaged vector field method
}

\author{Haochen Li         \and
        Yushun Wang       \\ \and
        Mengzhao Qin
}

\institute{H. Li \at
              Jiangsu Key Laboratory for NSLSCS, School of Mathematical Sciences, Nanjing Normal University, Jiangsu, $210023$, People's Republic of China
\and
Y. Wang \at
              Jiangsu Key Laboratory for NSLSCS, School of Mathematical Sciences, Nanjing Normal University, Jiangsu, $210023$, People's Republic of China\\
              Tel.: +86-13912961556\\
              Fax: +86-517-83525012\\
              \email{wangyushun@njnu.edu.cn}
           \and
           M. Qin  \at
              Lsec, Academy of Mathematics and System Sciences, Chinese Academy of Science, PoBox $2719$, Beijing, People's Republic of China
}

\maketitle

\begin{abstract}
In this paper, based on the theory of rooted trees and B-series,
we propose the concrete formulas of the substitution law for the trees of order $=5$.
With the help of the new substitution law, we derive a B-series integrator extending the averaged vector field (AVF) method to high order.
The new integrator turns out to be of order six and exactly preserves energy for Hamiltonian systems. Numerical experiments are presented to demonstrate the accuracy and the energy-preserving property of the sixth order AVF method.\\\\
\end{abstract}

\noindent
\textbf{Mathematics Subject Classfication} \quad 65D17 $\cdot$ 65L05 $\cdot$ 65P10

\section{Introduction}
\label{intro}
Geometric numerical integration methods have come to the fore,
partly as an alternative to traditional methods such as Runge-Kutta methods.
A numerical method is called geometric if it preserves one or more physical/geometric
properties of the system exactly (i.e. up to round-off error)~\cite{Ref3,Ref1,Ref2,Ref4}. Examples of such geometric
properties that can be preserved are (first) integrals, symplectic structures, symmetries and reversing symmetries, phase-space volumes, Lyapunov functions, foliations, etc.
Geometric methods have applications in many areas of physics, including celestial mechanics,
particle accelerators, molecular dynamics, fluid dynamics, pattern formation, plasma physics, reaction-diffusion equations, and meteorology~\cite{Ref9,Ref5,Ref6,Ref10,RefQH,Ref7,Ref8,RefW}.

In fact, the conservation of the energy function is one of the most relevant features characterizing a Hamiltonian system.
Methods that exactly preserve energy have been considered since several decades.
Many energy-preserving methods have been proposed~\cite{Ref14,Ref11,Ref13,Ref15,Ref12}.
The discrete gradient method is among the most popular methods for designing integral
preserving schemes for ordinary differential equations, which was perhaps first discussed
by O. Gonzalez~\cite{Ref16}. T. Matsuo proposed discrete variational method for nonlinear wave equation~\cite{Ref17}. L. Brugnano and F. Iavernaro proposed Hamiltonian boundary value methods~\cite{Ref18,Ref19}. More recently, the existence of energy-preserving B-series methods has been shown in \cite{Ref20}, and a practical integrator which is the
averaged vector field (AVF) method of order two has been proposed~\cite{Ref14,Ref11,Ref21,Ref22}. This method exactly preserves the energy of Hamiltonian systems, and in contrast to projection-type integrators, only requires evaluations of the vector field. It is symmetric and its Taylor
series has the structure of a B-series. For polynomial Hamiltonians, the integral can be evaluated
exactly, and the implementation is comparable to that of the implicit mid-point rule~\cite{Ref12}.

In recent years, there has been growing interest in high-order AVF methods, and the second, third and fourth order AVF methods have been proposed~\cite{Ref22}.
It is shown that the theory of B-series and the substitution law obtained by substituting a B-series into the vector field appearing in another B-series play an important role in constructing high order methods~\cite{Ref20,Ref5}. The substitution law for the trees of order $\leq 4$ has been shown in \cite{Ref23,Ref24,Ref5}.
The fourth order AVF method is obtained by the concrete formulas of the substitution law for the trees of order $\leq 4$.
To construct a B-series method which not only has high order accuracy but also preserves the Hamiltonian is an important and interesting topic.
However, the concrete formulas of the substitution law for the trees of order $\geq 5$ have not been proposed, as the corresponding calculations are sufficiently complicated.

There are two aims in this paper. The first aim is to propose the concrete formulas of the substitution law for the trees of order $=5$. As we know, a low order B-series integrator can be extend to high order by the substitution law.
Using the new obtained concrete formulas of the substitution law for the trees of order $=5$, one can extend a low order geometric B-series integrator to sixth order naturally, easily and automatically, such as the symplectic integrator, the energy-preserving integrator, the momentum-preserving integrator and so on. Using the new obtained substitution law, we also easily obtained the same sixth order symplectic integrator in~\cite{Ref24}.
The second aim is to derive a sixth order AVF method for Hamiltonian systems. By expanding the second order AVF method into a B-series and considering the substitution law for the trees of order $=5$, a new method can be constructed. The new method is derived by the concrete formulas of the substitution law for the trees of order $=5$.
We prove that the new method is of order six and can also preserve the energy of Hamiltonian systems exactly.

The paper is organized as follows: In Sect.~\ref{sec:2}, we introduce the AVF method.
In Sect.~\ref{sec:3}, we recall a few definitions and properties related to trees and B-series. The substitution law for the trees of order $\leq 4$ is shown and we obtain that for the trees of order $=5$.
In Sect.~\ref{sec:4}, we derive the sixth order AVF method, and prove that the new method is of order six and it can preserve the Hamiltonian.
A few numerical experiments are
given in Sect.~\ref{sec:5} to confirm the theoretical results.
We finish the paper with conclusions in Sect.~\ref{sec:6}.

\section{The AVF method and its energy conservation property}
\label{sec:2}

Here we briefly discuss the AVF method and its energy-conservation property.
We consider a Hamiltonian differential equation, written in the form
\begin{equation}\label{eq:2s1}
\dot{z}=f(z)=S\nabla H(z),\quad z(0)=z_{0},\quad z\in \mathbb{R}^{n},
\end{equation}
where $f(z):\mathbb{R}^{n}\rightarrow \mathbb{R}$,
$S$ is a skew-symmetric constant matrix, $n$ is an even number,
and the Hamiltonian $H(z)$ is assumed to be sufficiently differentiable.
From system \eqref{eq:2s1}, we can get
\begin{equation}\label{eq:2s2}
\frac{dH(z(t))}{dt}=\nabla H(z)^{T}f(z)=\nabla H(z)^{T}S\nabla H(z)=0.
\end{equation}
Therefore, the flow of the system \eqref{eq:2s1} preserves the Hamiltonian $H(z)$ exactly.

The so-called AVF method~\cite{Ref14,Ref22} is defined by
\begin{equation}\label{eq:2s3}
\frac{z_{1}-z_{0}}{h}=\int_{0}^{1}f(\xi z_{1}+(1-\xi)z_{0})d\xi,
\end{equation}
where $h$ is the time step.

The AVF method \eqref{eq:2s3} can be rewritten as
\begin{equation*}
\frac{z_{1}-z_{0}}{h}=S\int_{0}^{1}\nabla H(\xi z_{1}+(1-\xi)z_{0})d\xi.
\end{equation*}
We can obtain
\begin{eqnarray*}
& & \frac{1}{h}(H(z_{1})-H(z_{0}))=\frac{1}{h}\int_{0}^{1}\frac{d}{d\xi}H(\xi z_{1}+(1-\xi)z_{0})d\xi\\
& = & (\int_{0}^{1}\nabla H(\xi z_{1}+(1-\xi)z_{0})d\xi)^{T}(\frac{z_{1}-z_{0}}{h})\\
& = & (\int_{0}^{1}\nabla H(\xi z_{1}+(1-\xi)z_{0})d\xi)^{T} S\int_{0}^{1}\nabla H(\xi z_{1}+(1-\xi)z_{0})d\xi=0.
\end{eqnarray*}
It follows that the Hamiltonian $H$ is conserved at every time step.

R. I. McLachlan and G. R. W. Quispel also proposed a higher order energy-preserving method~\cite{Ref22}
\begin{equation}\label{eq:2s4}
\frac{z_{1}^{i}-z_{0}^{i}}{h}=
(\delta_{j}^{i}+\alpha h^{2}f_{k}^{i}(\hat{z})f_{j}^{k}(\hat{z}))\int_{0}^{1}f^{j}(\xi z_{1}+(1-\xi)z_{0})d\xi,
\end{equation}
where $\alpha$ is an arbitrary constant, $\delta_{i}^{j}$ is the Kronecker delta, and we can take e.g. $\hat{z}=z_{0}$ or $\hat{z}=\frac{z_{0}+z_{1}}{2}$.
For $\alpha=0$, we recover the second order method \eqref{eq:2s3}. For $\alpha=-\frac{1}{12}$ and $\hat{z}=z_{0}$, the method is of order three. For $\alpha=-\frac{1}{12}$ and $\hat{z}=\frac{z_{0}+z_{1}}{2}$, the method is of order four.

The highest order of existing AVF method is four, and based on the theory of substitution law we can obtain a new AVF method of order six.

\section{Preliminaries}
\label{sec:3}

In this section, we briefly recall a few definitions and properties related to rooted trees and B-series~\cite{Ref25,Ref23,Ref24,Ref27,Ref20,Ref5,Ref6,Ref2,Ref26}.

\subsection{Trees and B-series}

Let $\emptyset$ denote the empty tree.

\begin{definition}[Unordered trees~\cite{Ref5}] The set $\mathcal{T}$ of (rooted) unordered trees is recursively defined by
\begin{equation}\label{eq:3s1}
\tree11 \in \mathcal{T},\quad [\tau_{1},\ldots,\tau_{m}]\in \mathcal{T},\quad \forall~ \tau_{1},\ldots,\tau_{m}\in \mathcal{T},
\end{equation}
where $\tree11$ is the tree with only one vertex, and $\tau=[\tau_{1},\ldots,\tau_{m}]$ represents the rooted tree obtained by grafting the
roots of $\tau_{1},\ldots,\tau_{m}\in \mathcal{T}$ to a new vertex. Trees $\tau_{i}$ are called the branches of $\tau$.
\end{definition}
We note that $\tau$ does not depend on the ordering of $\tau_{1},\ldots,\tau_{m}$. For instance, we do not distinguish between $[\tree11,[\tree11]]$ and $[[\tree11],\tree11]$.

\begin{definition}[Coefficients~\cite{Ref2}] The following coefficients are defined recursively for all trees $\tau=[\tau_{1},\ldots,\tau_{m}]\in \mathcal{T}$:
\begin{align*}
|\tau|~=~ & 1+\sum_{i=1}^{m}|\tau_{i}|\quad & (the~order,~i.e.~the~number~of~vertices),\\
\alpha(\tau)~=~ & \frac{(|\tau|-1)!}{|\tau_{1}|!\cdot\ldots\cdot|\tau_{m}|!}\alpha(\tau_{1})\cdot\ldots\cdot\alpha(\tau_{m})\frac{1}{\mu_{1}!\mu_{2}!\ldots}\quad & (Connes-Moscovici~ weights),\\
\sigma(\tau)~=~ & \alpha(\tau_{1})\cdot\ldots\cdot\alpha(\tau_{m})\cdot\mu_{1}!\mu_{2}!\cdot\ldots\quad & (symmetry),\\
\gamma(\tau)~=~ & |\tau|\gamma(\tau_{1})\cdot\ldots\cdot\gamma(\tau_{m})\quad & (density),
\end{align*}
where the integers $\mu_{1},\mu_{2},\ldots$ count equal trees among $\tau_{1},\ldots,\tau_{m}$.
\end{definition}

\begin{definition}[Elementary differentials~\cite{Ref23}] For a vector field $f:\mathbb{R}^{d}\rightarrow\mathbb{R}^{d}$, and for an unordered tree
$\tau=[\tau_{1},\ldots,\tau_{m}]\in \mathcal{T}$, the so-called elementary differential is a mapping $F_{f}(\tau):\mathbb{R}^{d}\rightarrow\mathbb{R}^{d}$, recursively defined by
\begin{equation*}
F_{f}(\tree11)(z)=f(z),\quad F_{f}(\tau)(z)=f^{m}(z)(F_{f}(\tau_{1})(z),\ldots,F_{f}(\tau_{m})(z)).
\end{equation*}
\end{definition}

\begin{definition}[B-series~\cite{Ref26}] For a mapping $a:\mathcal{T}\bigcup\{\emptyset\}\rightarrow\mathbb{R}$, a formal series of the form
\begin{equation*}
B_{f}(a,z)=a(\emptyset)z+\sum_{\tau\in \mathcal{T}}\frac{h^{|\tau|}}{\sigma(\tau)}a(\tau)F_{f}(\tau)(z)
\end{equation*}
is called a B-series.
\end{definition}

\begin{theorem}[Exact solution~\cite{Ref2}]
If $z(t)$ denotes the exact solution of \eqref{eq:2s1}, it holds for all $j\geq 1$,
\begin{equation*}
\frac{1}{j!}z^{(j)}(0)=\sum_{\tau\in \mathcal{T}, |\tau|=j}\frac{1}{\sigma(\tau)\gamma(\tau)}F_{f}(\tau)(z_{0}).
\end{equation*}
Therefore, letting $\gamma(\emptyset)=1$, the exact solution of \eqref{eq:2s1} is (formally) given by
\begin{equation*}
z(h)=B_{f}(\frac{1}{\gamma},z_{0}).
\end{equation*}
\end{theorem}

\subsection{Basic tools for trees}

\subsubsection{Partitions and skeletons}

In order to manipulate trees more conveniently, it is useful to consider the set $\mathcal{OT}$ of ordered trees defined below.\\

\begin{definition}[Ordered trees~\cite{Ref2}]\label{def:3s5} The set $\mathcal{OT}$ of ordered trees is recursively defined by
\begin{equation*}
\tree11 \in \mathcal{OT},\quad (\omega_{1},\ldots,\omega_{m})\in \mathcal{OT},\quad \forall~ \omega_{1},\ldots,\omega_{m}\in \mathcal{OT}.
\end{equation*}
In contrast to $\mathcal{T}$, the ordered tree $\omega$ depends on the ordering $\omega_{1},\ldots,\omega_{m}$.
\end{definition}

Neglecting the ordering, a tree $\tau \in \mathcal{T}$ can be considered as an equivalent class of ordered trees, denoted $\tau=\overline{\omega}$.
Therefore, any function $\psi$ defined on $\mathcal{T}$ (such as order, symmetry, density,\ldots) can be extended to $\mathcal{OT}$ by
putting $\psi(\omega)=\psi(\overline{\omega})$ for all $\omega \in \mathcal{OT}$. Moreover, for all  $\tau \in \mathcal{T}$, we can choose a tree $\omega(\tau) \in \mathcal{OT}$ such as $\tau=\overline{\omega(\tau)}$~\cite{Ref23}.\\

\begin{definition}[Partitions of a tree~\cite{Ref23}] A partition $p^{\theta}$ of an ordered tree $\theta \in \mathcal{OT}$ is the (ordered) tree obtained
from $\theta$ by replacing some of its edges by dashed ones. We denote $P(p^{\theta}) = \{ s_{1},\ldots,s_{k}\}$ the list of subtrees $s_{i}\in \mathcal{T}$
obtained from $p^{\theta}$ by removing dashed edges and by neglecting the ordering of each subtree. We denote $\#(p^{\theta})=k$
the number of $s_{i}$'s. We observe that precisely one of the $s_{i}$'s contains the root of $\theta$. We denote this distinguished
tree by $R(p^{\theta})\in \mathcal{T}$. We denote $P^{\ast}(p^{\theta}) = P(p^{\theta})\backslash \{ R(p^{\theta})\}$ the list of $s_{i}$'s that do not contain the root of $\theta$. Eventually, the set of all partitions $p^{\theta}$ of $\theta$ is denoted $\mathcal{P}(\theta)$. Finally, for $\tau \in \mathcal{T}$, we put $\mathcal{P}(\tau)=\mathcal{P}(\omega(\tau))$ where $\omega(\tau) \in \mathcal{OT}$ is given in definition \ref{def:3s5}.
\end{definition}

We observe that any tree $\tau \in \mathcal{T}$ has exactly $2^{|\tau|-1}$ partitions $p^{\tau}\in\mathcal{P}(\tau)$, and that different partitions may lead
to the same list of subtrees $P(p^{\tau})$.\\

\begin{definition}[Skeleton of a partition~\cite{Ref27}]
The skeleton $\chi(p^{\tau}) \in \mathcal{T}$ of a partition $p^{\tau}\in\mathcal{P}(\tau)$ of a tree $\tau \in \mathcal{T}$ is the tree obtained by replacing in $p^{\tau}$ each tree of $P(p^{\tau})$ by a single vertex and then dashed edges by solid ones. We can notice that $|\chi(p^{\tau})|=\#(p^{\tau})$.
\end{definition}

\subsubsection{Substitution law}

\begin{theorem}[Substitution law~\cite{Ref23,Ref24,Ref5}]
Let $a,b:\mathcal{T}\bigcup\{\emptyset\}\rightarrow\mathbb{R}$ be two mappings with $b(\emptyset) = 0$. Given a field $f:\mathbb{R}^{d}\rightarrow\mathbb{R}^{d}$, consider the (h-dependent)
field $g:\mathbb{R}^{d}\rightarrow\mathbb{R}^{d}$ defined by
\begin{equation*}
hg(z)=B_{f}(b,z),
\end{equation*}
Then, there exists a mapping $b\star a:\mathcal{T}\bigcup\{\emptyset\}\rightarrow\mathbb{R}$ satisfying
\begin{equation*}
B_{g}(a,z)=B_{f}(b\star a,z),
\end{equation*}
and $b\star a$ is defined by
\begin{equation}\label{eq:3s2}
b\star a(\emptyset)=a(\emptyset),\quad \forall \tau \in \mathcal{T},b\star a(\tau)=\sum_{p^{\tau}\in\mathcal{P}(\tau)}a(\chi(p^{\tau}))\prod_{\delta\in P(p^{\tau})}b(\delta).
\end{equation}
\end{theorem}

in \cite{Ref23,Ref5}, the concrete formulas of the substitution law for the trees of order $\leq 4$ was proposed (see table~\ref{tab:1}).
In this paper, we obtain the concrete formulas of the substitution law for the trees of order $=5$ (see table~\ref{tab:2}).

\begin{table}
\caption{The concrete formulas of the substitution law $\star$ defined in \eqref{eq:3s2} for trees of order $\leq 4$}
\label{tab:1}
\begin{align*}
b\star a(\emptyset) =~ & a(\emptyset),\\
b\star a(\tree11) =~ & a(\tree11)b(\tree11),  \\
b\star a(\tree21) =~ & a(\tree11)b(\tree21) + a(\tree21)b(\tree11)^{2},  \\
b\star a(\tree31) =~ & a(\tree11)b(\tree31) + 2a(\tree21)b(\tree11)b(\tree21) + a(\tree31)b(\tree11)^{3},  \\\\
b\star a(\tree32) =~ & a(\tree11)b(\tree32) + 2a(\tree21)b(\tree11)b(\tree21) + a(\tree32)b(\tree11)^{3},  \\
b\star a(\tree41) =~ & a(\tree11)b(\tree41) + 3a(\tree21)b(\tree11)b(\tree31) + 3a(\tree31)b(\tree11)^{2}b(\tree21) + a(\tree41)b(\tree11)^{4},  \\\\
b\star a(\tree42) =~ & a(\tree11)b(\tree42) + a(\tree21)b(\tree11)b(\tree32) + a(\tree21)b(\tree21)^{2} + a(\tree21)b(\tree11)b(\tree31) \\\\
                    & + 2a(\tree31)b(\tree11)^{2}b(\tree21) + a(\tree32)b(\tree11^{2})b(\tree21) + a(\tree42)b(\tree11)^{4},  \\\\
b\star a(\tree43) =~ & a(\tree11)b(\tree43) + a(\tree21)b(\tree11)b(\tree31) + 2a(\tree21)b(\tree11)b(\tree32) + a(\tree31)b(\tree11)^{2}b(\tree21) \\\\
                    & + 2a(\tree32)b(\tree11)^{2}b(\tree21) + a(\tree43)b(\tree11)^{4},  \\\\
b\star a(\tree44) =~ & a(\tree11)b(\tree44) + 2a(\tree21)b(\tree11)b(\tree32) + a(\tree21)b(\tree21)^{2} + 3a(\tree32)b(\tree11)^{2}b(\tree21)\\\\
                    & + a(\tree44)b(\tree11)^{4}.
\end{align*}
\end{table}

\begin{table}
\caption{The concrete formulas of the substitution law $\star$ defined in \eqref{eq:3s2} for trees of order $=5$}
\label{tab:2}
\begin{align*}
b\star a(\tree51) =~ & a(\tree11)b(\tree51) + 4a(\tree21)b(\tree11)b(\tree41) + 6a(\tree31)b(\tree11)^{2}b(\tree31) \\
                    & + 4a(\tree41)b(\tree11)^{3}b(\tree21) + a(\tree51)b(\tree11)^{5},\\\\
b\star a(\tree52) =~ & a(\tree11)b(\tree52) + a(\tree21)b(\tree11)b(\tree41) + 2a(\tree21)b(\tree11)b(\tree42) + a(\tree21)b(\tree21)b(\tree31) \\\\
                     & + a(\tree32)b(\tree11)^{2}b(\tree31) + 2a(\tree31)b(\tree11)^{2}b(\tree31) + a(\tree31)b(\tree11)^{2}b(\tree32) \\\\
                     & + 2a(\tree31)b(\tree11)b(\tree21)^{2} + 2a(\tree41)b(\tree11)^{3}b(\tree21) + 2a(\tree42)b(\tree11)^{3}b(\tree21) \\\\
                     & + a(\tree52)b(\tree11)^{5},\\\\
b\star a(\tree53) =~ & a(\tree11)b(\tree53) + 2a(\tree21)b(\tree11)b(\tree42) + 2a(\tree21)b(\tree21)b(\tree32) \\\\
                     & + a(\tree31)b(\tree11)^{2}b(\tree31) + 2a(\tree32)b(\tree11)^{2}b(\tree32) + 3a(\tree31)b(\tree11)b(\tree21)^{2} \\\\
                     & + 4a(\tree42)b(\tree11)^{3}b(\tree21) + a(\tree53)b(\tree11)^{5},\\\\
b\star a(\tree54) =~ & a(\tree11)b(\tree54) + a(\tree21)b(\tree11)b(\tree41) + 3a(\tree21)b(\tree11)b(\tree43) \\\\
                     & + 3a(\tree32)b(\tree11)^{2}b(\tree31) + 3a(\tree31)b(\tree11)^{2}b(\tree32) + a(\tree41)b(\tree11)^{3}b(\tree21) \\\\
                     & + 3a(\tree43)b(\tree11)^{3}b(\tree21) + a(\tree54)b(\tree11)^{5},\\\\
b\star a(\tree55) =~ & a(\tree11)b(\tree55) + a(\tree21)b(\tree11)b(\tree42) + a(\tree21)b(\tree11)b(\tree44) + a(\tree21)b(\tree11)b(\tree43) \\\\
                     & + a(\tree21)b(\tree21)b(\tree32) + 2a(\tree32)b(\tree11)^{2}b(\tree32) + a(\tree31)b(\tree11)b(\tree21)^{2} \\\\
                     & + a(\tree32)b(\tree11)b(\tree21)^{2} + a(\tree32)b(\tree11)^{2}b(\tree31) + a(\tree31)b(\tree11)^{2}b(\tree32) \\\\
                     & + a(\tree42)b(\tree11)^{3}b(\tree21) + a(\tree44)b(\tree11)^{3}b(\tree21) + 2a(\tree43)b(\tree11)^{3}b(\tree21) \\\\
                     & + a(\tree55)b(\tree11)^{5},\\\\
b\star a(\tree56) =~ & a(\tree11)b(\tree56) + a(\tree21)b(\tree11)b(\tree42) + a(\tree21)b(\tree21)b(\tree31) + a(\tree21)b(\tree21)b(\tree32) \\\\
                     & + a(\tree21)b(\tree11)b(\tree44) + a(\tree32)b(\tree11)^{2}b(\tree31) + 2a(\tree32)b(\tree11)b(\tree21)^{2} \\\\
                     & + a(\tree31)b(\tree11)b(\tree21)^{2} + 2a(\tree31)b(\tree11)^{2}b(\tree32) + a(\tree44)b(\tree11)^{3}b(\tree21) \\\\
                     & + 3a(\tree42)b(\tree11)^{3}b(\tree21) + a(\tree56)b(\tree11)^{5},\\
\end{align*}
\end{table}
\begin{table}
\begin{align*}
b\star a(\tree57) =~ & a(\tree11)b(\tree57) + a(\tree21)b(\tree11)b(\tree43) + 2a(\tree21)b(\tree11)b(\tree44) + a(\tree21)b(\tree21)b(\tree31) \\\\
                     & + 2a(\tree32)b(\tree11)b(\tree21)^{2} + 2a(\tree32)b(\tree11)^{2}b(\tree32) + a(\tree32)b(\tree11)^{2}b(\tree31) \\\\
                     & + a(\tree31)b(\tree11)^{2}b(\tree32) + 2a(\tree43)b(\tree11)^{3}b(\tree21) + 2a(\tree44)b(\tree11)^{3}b(\tree21) \\\\
                     & + a(\tree57)b(\tree11)^{5},\\\\
b\star a(\tree58) =~ & a(\tree11)b(\tree58) + 2a(\tree21)b(\tree11)b(\tree44) + 2a(\tree21)b(\tree21)b(\tree32) + 3a(\tree32)b(\tree11)^{2}b(\tree32) \\\\
                     & + 3a(\tree32)b(\tree11)b(\tree21)^{2} + 4a(\tree44)b(\tree11)^{3}b(\tree21) + a(\tree58)b(\tree11)^{5} \\\\
b\star a(\tree59) =~ & a(\tree11)b(\tree59) + 2a(\tree21)b(\tree11)b(\tree42) + a(\tree21)b(\tree21)b(\tree31) + a(\tree21)b(\tree11)b(\tree43) \\\\
                     & + 2a(\tree31)b(\tree11)^{2}b(\tree31) + 2a(\tree31)b(\tree11)^{2}b(\tree32) + 2a(\tree32)b(\tree11)b(\tree21)^{2} \\\\
                     & + a(\tree43)b(\tree11)^{3}b(\tree21) + 2a(\tree42)b(\tree11)^{3}b(\tree21) + a(\tree41)b(\tree11)^{3}b(\tree21) \\\\
                     & + a(\tree59)b(\tree11)^{5}.
\end{align*}
\end{table}

\section{Sixth order AVF method}
\label{sec:4}
\subsection{The second order AVF method and its B-series}

Consider an ordinary differential equation
\begin{equation}\label{eq:4s1}
\dot{z}=f(z),\quad z\in\mathbb{R}^{n},\quad z(t_{0})=z_{0},
\end{equation}
and the second order AVF method
\begin{equation}\label{eq:4s2}
\Phi_{h}^{f}(z_{0})=z_{1}=z_{0}+h\int_{0}^{1}f(\xi z_{1}+(1-\xi)z_{0})d\xi.
\end{equation}
\begin{theorem} The AVF method \eqref{eq:4s2} can be expanded into a B-series
\begin{equation}\label{eq:4s3}
\Phi_{h}^{f}(z_{0})=B_{f}(a,z_{0})=a(\emptyset)z_{0}+\sum_{\tau\in \mathcal{T}}\frac{h^{|\tau|}}{\sigma(\tau)}a(\tau)F_{f}(\tau)(z_{0}),
\end{equation}
where $a(\emptyset)=a(\tree11)=1$, and for all $\tau=[\tau_{1},\ldots,\tau_{m}]\in \mathcal{T}$,
\begin{equation*}
a(\tau)=\frac{1}{m+1}a(\tau_{1})\cdots a(\tau_{m}).
\end{equation*}
\end{theorem}
\paragraph{Proof} We develop the derivatives of \eqref{eq:4s2}, by Leibniz's rule, and obtain
\begin{equation*}
z_{1}^{(q)}=[h\int_{0}^{1}f(\xi z_{1}+(1-\xi)z_{0})d\xi]^{(q)}=h[\int_{0}^{1}f(\xi z_{1}+(1-\xi)z_{0})d\xi]^{(q)}+q[\int_{0}^{1}f(\xi z_{1}+(1-\xi)z_{0})d\xi]^{q-1}.
\end{equation*}
This gives, for $h=0$,
\begin{equation*}
z^{(q)}:=z_{1}^{(q)}|_{h=0}=q[\int_{0}^{1}f(\xi z_{1}+(1-\xi)z_{0})d\xi]^{q-1}|_{h=0},\quad q\geq 1,
\end{equation*}
and considering $z_{1}|_{h=0}=z_{0}$, we can obtain
\begin{align}\label{eq:4s4}
\dot{z}=~ & \int_{0}^{1}f(\xi z_{0}+(1-\xi)z_{0})d\xi=1\cdot 1\cdot 1\cdot f(z_{0}),\nonumber\\
\ddot{z}=~ & 2\int_{0}^{1}\xi f'(\xi z_{0}+(1-\xi)z_{0})\dot{z}d\xi=f(z_{0})=2\cdot 1\cdot \frac{1}{2}f'(z_{0})\dot{z},\nonumber\\
z^{(3)}=~ & 3\int_{0}^{1}(\xi^{2} f''(\xi z_{0}+(1-\xi)z_{0})(\dot{z},\dot{z})+\xi f'(\xi z_{0}+(1-\xi)z_{0})\ddot{z})d\xi \nonumber\\
       =~ & 3\cdot(1\cdot \frac{1}{3}f''(z_{0})(\dot{z},\dot{z})+1\cdot \frac{1}{2}f'(z_{0})\ddot{z}),\\
z^{(4)}=~ & 4\int_{0}^{1}(\xi^{3} f'''(\xi z_{0}+(1-\xi)z_{0})(\dot{z},\dot{z},\dot{z})+3\xi^{2} f''(\xi z_{0}+(1-\xi)z_{0})(\ddot{z},\dot{z})\nonumber\\
          & +\xi f'(\xi z_{0}+(1-\xi)z_{0})z^{(3)})d\xi \nonumber\\
       =~ & 4\cdot(1\cdot \frac{1}{4}f'''(z_{0})(\dot{z},\dot{z},\dot{z})+3\cdot \frac{1}{3}f''(z_{0})(\ddot{z},\dot{z})+1\cdot \frac{1}{2}f'(z_{0})z^{(3)}),\nonumber
\end{align}
and so on. We now insert in \eqref{eq:4s4} recursively the computed derivatives $\dot{z},\ddot{z},\ldots$ into the right side of the subsequent formulas.
Letting $f^{(q)}:=f^{(q)}(z_{0})$ and $F(\tau)=F_{f}(\tau)(z_{0})$, we can obtain
\begin{align}\label{eq:4s5}
\dot{z}=~ & 1\cdot 1\cdot 1\cdot f=\gamma(\tree11)\alpha(\tree11)a(\tree11)F(\tree11),\nonumber\\
\ddot{z}=~ & 2\cdot 1\cdot \frac{1}{2}f'f=\gamma(\tree21)\alpha(\tree21)a(\tree21)F(\tree21),\nonumber\\
z^{(3)}=~ & 3\cdot 1\cdot \frac{1}{3}f''(f,f)+(2\cdot 3)\cdot 1\cdot (\frac{1}{2}\cdot \frac{1}{2})f'f'f \nonumber\\
       ~ & =\gamma(\tree31)\alpha(\tree31)a(\tree31)F(\tree31)+\gamma(\tree32)\alpha(\tree32)a(\tree32)F(\tree32),\\
z^{(4)}=~ & 4\cdot 1\cdot \frac{1}{4}f'''(f,f,f)+(2\cdot 4)\cdot 3\cdot (\frac{1}{2}\cdot \frac{1}{3})f''(f'f,f)+(3\cdot 4)\cdot 1\cdot (\frac{1}{2}\cdot \frac{1}{3})f'f''(f,f)\nonumber\\
          & +(2\cdot 3\cdot 4)\cdot 1\cdot (\frac{1}{2}\cdot \frac{1}{2}\cdot \frac{1}{2})f'f'f'f\nonumber\\
       =~ & \gamma(\tree41)\alpha(\tree41)a(\tree41)F(\tree41)+\gamma(\tree42)\alpha(\tree42)a(\tree42)F(\tree42)\nonumber\\ \nonumber\\
          & +\gamma(\tree43)\alpha(\tree43)a(\tree43)F(\tree43)+\gamma(\tree44)\alpha(\tree44)a(\tree44)F(\tree44),\nonumber
\end{align}
and so on. For all $\tau=[\tau_{1},\ldots,\tau_{m}]\in \mathcal{T}$, letting $\alpha=\frac{(|\tau|-1)!}{|\tau_{1}|!\cdot\ldots\cdot|\tau_{m}|!}\frac{1}{\mu_{1}!\mu_{2}!\ldots}$,
where the integers $\mu_{1},\mu_{2},\ldots$ count equal trees among $\tau_{1},\ldots,\tau_{m}$, we can obtain
\begin{align*}
& \gamma(\tau)\alpha(\tau)a(\tau)F(\tau)\\
=~ & |\tau|\int_{0}^{1}\alpha\xi^{m} f^{m}(z_{0})(\gamma(\tau_{1})\alpha(\tau_{1})a(\tau_{1})F(\tau_{1}),\ldots,\gamma(\tau_{m})\alpha(\tau_{m})a(\tau_{m})F(\tau_{m}))d\xi\\
=~ & [|\tau|\gamma(\tau_{1})\cdot\ldots\cdot\gamma(\tau_{m})]\cdot [\alpha\alpha(\tau_{1})\cdot\ldots\cdot\alpha(\tau_{m})]\cdot [\frac{1}{m+1}a(\tau_{1})\cdots a(\tau_{m})]\cdot f^{m}(z_{0})(F(\tau_{1}),\ldots,F(\tau_{m}))\\
=~ & \gamma(\tau)\alpha(\tau)[\frac{1}{m+1}a(\tau_{1})\cdots a(\tau_{m})]F(\tau).
\end{align*}
So we have
\begin{equation*}
z^{(q)}=\sum_{\tau\in \mathcal{T},|\tau|=q}\gamma(\tau)\alpha(\tau)a(\tau)F(\tau),
\end{equation*}
where $a(\tree11)=1$, and for all $\tau=[\tau_{1},\ldots,\tau_{m}]\in \mathcal{T}$,
\begin{equation*}
a(\tau)=\frac{1}{m+1}a(\tau_{1})\cdots a(\tau_{m}).
\end{equation*}
Letting $a(\emptyset)=1$, and considering $\sigma(\tau)=\frac{|\tau|!}{\alpha(\tau)\gamma(\tau)}$, we obtain
\begin{equation*}
\Phi_{h}^{f}(z_{0})=z_{1}=a(\emptyset)z_{0}+\sum_{\tau\in \mathcal{T}}\frac{h^{|\tau|}}{\sigma(\tau)}a(\tau)F_{f}(\tau)=B_{f}(a,z_{0}).
\end{equation*}
The proof is completed.

The B-series \eqref{eq:4s3} can be rewritten as
\begin{align*}
\Phi_{h}^{f}(z_{0})=~ & z_{0}+hf+\frac{1}{2}h^{2}f'f+h^{3}(\frac{1}{2\cdot 3}f''(f,f)+\frac{1}{4}f'f'f)+h^{4}(\frac{1}{4\cdot 6}f'''(f,f,f)\\
                    & +\frac{1}{6}f''(f'f,f)+\frac{1}{2\cdot 6}f'f''(f,f)+\frac{1}{8}f'f'f'f)+\ldots\\
                   =~ & z_{0}+hF(\tree11)+\frac{1}{2}h^{2}F(\tree21)+h^{3}(\frac{1}{2\cdot 3}F(\tree31)+\frac{1}{4}F(\tree32))+h^{4}(\frac{1}{4\cdot 6}F(\tree41) +\frac{1}{6}F(\tree42)\\
                     & +\frac{1}{2\cdot 6}F(\tree43)+\frac{1}{8}F(\tree44))+h^{5}(\frac{1}{5\cdot 24}F(\tree51) +\frac{1}{2\cdot 8}F(\tree52)+\frac{1}{2\cdot 12}F(\tree53)\\\\
                     & +\frac{1}{6\cdot 8}F(\tree54)+\frac{1}{12}F(\tree55)+\frac{1}{12}F(\tree56)+\frac{1}{2\cdot 12}F(\tree57)+\frac{1}{16}F(\tree58)+\frac{1}{2\cdot 9}F(\tree59))+\ldots
\end{align*}

\subsection{Sixth order AVF method}

Let $a:\mathcal{T}\bigcup\{\emptyset\}\rightarrow\mathbb{R}$ be a mapping satisfying $a(\emptyset)=1$, $a(\tree11)\neq 0$, and let $f:\mathbb{R}^{d}\rightarrow\mathbb{R}^{d}$ be a field. We consider here the numerical flow $\Phi_{h}^{g}$, where $g:\mathbb{R}^{d}\rightarrow\mathbb{R}^{d}$ denotes the modified field of $f$, whose B-series expansion is
\begin{equation*}
\Phi_{h}^{g}(z)=B_{g}(a,z).
\end{equation*}
The fundamental idea of obtaining sixth order AVF method consists in interpreting the numerical solution $z_{1}=\Phi_{h}^{g}$ of the
initial value problem $z(0)=z_{0}$, $\dot{z}=g(z)$ as the exact solution of a modified differential equation $\dot{z}=f(z)$~\cite{Ref20,Ref5}.\\

\begin{theorem} If $z_{1}=\Phi_{h}^{f}(z_{0})$, which is a numerical solution of \eqref{eq:4s1}, can be expanded into a B-series
\begin{equation*}
\Phi_{h}^{f}(z_{0})=z_{1}=B_{f}(a,z_{0}), a(\emptyset)=1,
\end{equation*}
and $b:\mathcal{T}\bigcup\{\emptyset\}\rightarrow\mathbb{R}$ be a mapping with $b(\emptyset) = 0$ satisfying
\begin{equation}\label{eq:4s6}
\forall \tau\in \mathcal{T}, b\star a(\tau)=\frac{1}{\gamma(\tau)},
\end{equation}
then we can obtain a new numerical solution $\Phi_{h}^{g}(z_{0})=B_{g}(a,z_{0})$, satisfying
\begin{equation*}
\Phi_{h}^{g}(z_{0})=\varphi_{h}^{f}(z_{0}),
\end{equation*}
where $\varphi_{h}^{f}(z_{0})$ denotes the exact solution of \eqref{eq:4s1}, and $g:\mathbb{R}^{n}\rightarrow\mathbb{R}^{n}$ is a (h-dependent) field defined by
\begin{equation*}
hg(y)=B_{f}(b,y).
\end{equation*}
\end{theorem}
\begin{proof} From theorem 3.1, $\varphi_{h}^{f}(z_{0})=B_{f}(\frac{1}{\gamma},z_{0})$.
And from theorem 3.2, we obtain
\begin{equation*}
\Phi_{h}^{g}(z_{0})=B_{g}(a,z_{0})=B_{f}(b\star a,z_{0})=B_{f}(\frac{1}{\gamma},z_{0})=\varphi_{h}^{f}(z_{0}).
\end{equation*}
The proof is completed.
\end{proof}

Letting $e(\tau)=\frac{1}{\gamma(\tau)}$, we can calculate $b(\tau)$ defined in \eqref{eq:4s6} for trees of order $\leq 6$:
\begin{align*}
& |\tau|=1, & \\
& a(\tree11)b(\tree11) =~ e(\tree11),\quad a(\tree11)=1, e(\tree11)=1,  &  b(\tree11)=1, \\
& |\tau|=2, & \\
& a(\tree11)b(\tree21) + a(\tree21)b(\tree11)^{2} =~ e(\tree21),\\
& a(\tree21)=\frac{1}{2}, e(\tree21)=\frac{1}{2},  &  b(\tree21)=0,\\
& |\tau|=3, & \\
& a(\tree11)b(\tree31) + 2a(\tree21)b(\tree11)b(\tree21) + a(\tree31)b(\tree11)^{3}=~ e(\tree31),\\
& b(\tree31) + \frac{1}{3}= \frac{1}{3},  &  b(\tree31)=0,\\
& a(\tree11)b(\tree32) + 2a(\tree21)b(\tree11)b(\tree21) + a(\tree32)b(\tree11)^{3}=~ e(\tree32),\\
& b(\tree32) + \frac{1}{4}= \frac{1}{6},  &  b(\tree31)=-\frac{1}{12},\\
& |\tau|=4, & \\
& a(\tree11)b(\tree41) + 3a(\tree21)b(\tree11)b(\tree31) + 3a(\tree31)b(\tree11)^{2}b(\tree21) + a(\tree41)b(\tree11)^{4}=~ e(\tree41),\\
& b(\tree41) + \frac{1}{4}= \frac{1}{4},  &  b(\tree31)=0,\\
& a(\tree11)b(\tree42) + a(\tree21)b(\tree11)b(\tree32) + a(\tree42)b(\tree11)^{4}=~ e(\tree42),\\
& b(\tree42) - \frac{1}{24}+ \frac{1}{6}= \frac{1}{8},  &  b(\tree42)=0,\\
& a(\tree11)b(\tree43) + 2a(\tree21)b(\tree11)b(\tree32) + a(\tree43)b(\tree11)^{4}=~ e(\tree43),\\
& b(\tree43) - \frac{1}{12}+ \frac{1}{6}= \frac{1}{12},  &  b(\tree43)=0,\\\\
& a(\tree11)b(\tree44) + 2a(\tree21)b(\tree11)b(\tree32)+ a(\tree44)b(\tree11)^{4}=~ e(\tree44),\\
& b(\tree44) - \frac{1}{12}+ \frac{1}{8}= \frac{1}{24},  &  b(\tree44)=0,
\end{align*}
\begin{align*}
& |\tau|=5, & \\
& a(\tree11)b(\tree51) + a(\tree51)b(\tree11)^{5}=~ e(\tree51),\\
& b(\tree51) + \frac{1}{5}= \frac{1}{5},  &  b(\tree51)=0,\\
& a(\tree11)b(\tree52) + a(\tree21)b(\tree11)b(\tree41) + a(\tree31)b(\tree11)^{2}b(\tree32) + a(\tree52)b(\tree11)^{5}=~ e(\tree52),\\
& b(\tree52) - \frac{1}{36} + \frac{1}{8}= \frac{1}{10},  &  b(\tree52)=\frac{1}{360},\\\\
& a(\tree11)b(\tree53) + 2a(\tree21)b(\tree21)b(\tree32) + a(\tree53)b(\tree11)^{5}=~ e(\tree53),\\
& b(\tree53) - \frac{1}{24} + \frac{1}{12}= \frac{1}{20},  &  b(\tree53)=\frac{1}{120},\\\\
& a(\tree11)b(\tree54) + 3a(\tree31)b(\tree11)^{2}b(\tree32) + a(\tree54)b(\tree11)^{5}=~ e(\tree54),\\
& b(\tree54) - \frac{1}{12} + \frac{1}{8}= \frac{1}{20},  &  b(\tree54)=\frac{1}{120},\\\\
& a(\tree11)b(\tree55) + 2a(\tree32)b(\tree11)^{2}b(\tree32) + a(\tree31)b(\tree11)^{2}b(\tree32) + a(\tree55)b(\tree11)^{5}=~ e(\tree55),\\\\
& b(\tree55) - \frac{1}{24} - \frac{1}{36} + \frac{1}{12}= \frac{1}{40},  &  b(\tree55)=\frac{1}{90},\\\\
& a(\tree11)b(\tree56) + 2a(\tree31)b(\tree11)^{2}b(\tree32) + a(\tree56)b(\tree11)^{5}=~ e(\tree56),\\
& b(\tree56) - \frac{1}{18} + \frac{1}{12}= \frac{1}{30},  &  b(\tree56)=\frac{1}{180},\\\\
& a(\tree11)b(\tree57) + 2a(\tree32)b(\tree11)^{2}b(\tree32) + a(\tree31)b(\tree11)^{2}b(\tree32) + a(\tree57)b(\tree11)^{5}=~ e(\tree57),\\\\
& b(\tree57) - \frac{1}{24} - \frac{1}{36} + \frac{1}{12}= \frac{1}{60},  &  b(\tree57)=\frac{1}{360},\\\\
& a(\tree11)b(\tree58) + 3a(\tree32)b(\tree11)^{2}b(\tree32) + a(\tree58)b(\tree11)^{5}=~ e(\tree58),\\\\
& b(\tree58) - \frac{1}{16} + \frac{1}{16}= \frac{1}{120},  &  b(\tree58)=\frac{1}{120},\\\\
& a(\tree11)b(\tree59) + 2a(\tree31)b(\tree11)^{2}b(\tree32) + a(\tree59)b(\tree11)^{5}=~ e(\tree59),\\
& b(\tree59) - \frac{1}{18} + \frac{1}{9}= \frac{1}{15},  &  b(\tree59)=\frac{1}{90},
\end{align*}
and in the same way, for all trees of order $=6$, we can obtain $b(\tau)=0$.

Then we obtain the exact modified field $g(z_{0})=\frac{1}{h}B_{f}(b,z_{0})$ and the sixth order modified field
\begin{align*}
\breve{g}(z_{0})=~ & f(z_{0})-\frac{h^{2}}{12}F(\tree32)+\frac{h^{4}}{720}[6F(\tree58)+4F(\tree59)+F(\tree57)+4F(\tree56)+F(\tree54)\\\\
                   & +F(\tree52)+3F(\tree53)+8F(\tree55)],
\end{align*}
satisfying $g(z_{0})=\breve{g}(z_{0})+\mathcal{O}(h^{6})$ and $\Phi_{h}^{g}(z_{0})=\varphi_{h}^{f}(z_{0})$.
And considering $\Phi_{h}^{\check{g}}(z_{0})=z_{1}=z_{0}+h\int_{0}^{1}\check{g}(\xi z_{1}+(1-\xi)z_{0})d\xi$, we have
\begin{equation*}
z_{1}-z_{0}=h\int_{0}^{1}f(\xi z_{1}+(1-\xi)z_{0})d\xi + \mathcal{O}(h^{3}).
\end{equation*}
Letting $\hat{z}:=\frac{z_{1}+z_{0}}{2}$, $f^{(q)}:=f^{(q)}(\hat{z})$, $F:=\int_{0}^{1}f(\xi z_{1}+(1-\xi)z_{0})d\xi$ and $\eta=\xi-\frac{1}{2}$,
we have $F=\int_{-\frac{1}{2}}^{\frac{1}{2}}f(\hat{z}+\eta(z_{1}-z_{0}))d\eta$ and for all $\tau \in \mathcal{T}$
\begin{equation*}
\int_{0}^{1}F_{f}(\tau)(\xi z_{1}+(1-\xi)z_{0})d\xi=\int_{-\frac{1}{2}}^{\frac{1}{2}}F_{f}(\tau)(\hat{z}+\eta(z_{1}-z_{0}))d\eta.
\end{equation*}
Considering $|\tau|=1$, we have $\tau=\tree11$ and
\begin{equation*}
\int_{0}^{1}F_{f}(\tau)(\xi z_{1}+(1-\xi)z_{0})d\xi=\int_{0}^{1}f(\xi z_{1}+(1-\xi)z_{0})d\xi=F.
\end{equation*}
We consider $\int_{-\frac{1}{2}}^{\frac{1}{2}}\eta d\eta=\int_{-\frac{1}{2}}^{\frac{1}{2}}\eta^{3} d\eta=0$, $\int_{-\frac{1}{2}}^{\frac{1}{2}}\eta^{2} d\eta=\frac{1}{12}$. For $|\tau|=3$, we have $\tau=\tree32$ and
\begin{align}\label{eq:4s7}
z_{1}-z_{0}=~ & h\int_{-\frac{1}{2}}^{\frac{1}{2}}f(\hat{z}+\eta(z_{1}-z_{0}))d\eta+\mathcal{O}(h^{3})=~ hF+\mathcal{O}(h^{3})\\ \label{eq:4s8}
or\quad z_{1}-z_{0} =~ & h\int_{-\frac{1}{2}}^{\frac{1}{2}}[f(\hat{z})+\eta f'(\hat{z})(z_{1}-z_{0})+\mathcal{O}(h^{2})]d\eta+\mathcal{O}(h^{3})\\ \label{eq:4s9}
=~ & hf(\hat{z})+\mathcal{O}(h^{3})=~ hf+\mathcal{O}(h^{3})\nonumber\\
or\quad z_{1}-z_{0} =~ & \mathcal{O}(h).
\end{align}
So we can obtain
\begin{align*}
& \int_{-\frac{1}{2}}^{\frac{1}{2}}F_{f}(\tree32)(\hat{z}+\eta(z_{1}-z_{0}))d\eta\\
=~ & \int_{-\frac{1}{2}}^{\frac{1}{2}}f'(\hat{z}+\eta(z_{1}-z_{0}))f'(\hat{z}+\eta(z_{1}-z_{0}))f(\hat{z}+\eta(z_{1}-z_{0}))d\eta\\
=~ & \int_{-\frac{1}{2}}^{\frac{1}{2}}[f'(\hat{z})+\eta f''(\hat{z})(z_{1}-z_{0})+\frac{\eta^{2}}{2!} f'''(\hat{z})(z_{1}-z_{0},z_{1}-z_{0})\\
& +\frac{\eta^{3}}{3!} f^{(4)}(\hat{z})(z_{1}-z_{0},z_{1}-z_{0},z_{1}-z_{0})+\mathcal{O}(h^{4})]\cdot [f'(\hat{z})+\eta f''(\hat{z})(z_{1}-z_{0})\\
& +\frac{\eta^{2}}{2!} f'''(\hat{z})(z_{1}-z_{0},z_{1}-z_{0})+\frac{\eta^{3}}{3!}f^{(4)}(\hat{z})(z_{1}-z_{0},z_{1}-z_{0},z_{1}-z_{0})+\mathcal{O}(h^{4})]\\
& \cdot f(\hat{z}+\eta(z_{1}-z_{0}))d\eta\\
=~ & \int_{-\frac{1}{2}}^{\frac{1}{2}}\{f'(\hat{z})f'(\hat{z})f(\hat{z}+\eta(z_{1}-z_{0}))+f'(\hat{z})[\eta f''(\hat{z})(z_{1}-z_{0})+\frac{\eta^{2}}{2!} f'''(\hat{z})(z_{1}-z_{0},z_{1}-z_{0})\\
& +\frac{\eta^{3}}{3!} f^{(4)}(\hat{z})(z_{1}-z_{0},z_{1}-z_{0},z_{1}-z_{0})]\cdot [f(\hat{z})+\eta f'(\hat{z})(z_{1}-z_{0})\\
& +\frac{\eta^{2}}{2!} f''(\hat{z})(z_{1}-z_{0},z_{1}-z_{0})+\mathcal{O}(h^{3})]+[\eta f''(\hat{z})(z_{1}-z_{0})+\frac{\eta^{2}}{2!} f'''(\hat{z})(z_{1}-z_{0},z_{1}-z_{0})\\
& +\frac{\eta^{3}}{3!} f^{(4)}(\hat{z})(z_{1}-z_{0},z_{1}-z_{0},z_{1}-z_{0})]\cdot [f'(\hat{z})+\eta f''(\hat{z})(z_{1}-z_{0})\\
& +\frac{\eta^{2}}{2!} f'''(\hat{z})(z_{1}-z_{0},z_{1}-z_{0})+\frac{\eta^{3}}{3!}f^{(4)}(\hat{z})(z_{1}-z_{0},z_{1}-z_{0},z_{1}-z_{0})]\\
& \cdot [f(\hat{z})+\eta f'(\hat{z})(z_{1}-z_{0})+\frac{\eta^{2}}{2!} f''(\hat{z})(z_{1}-z_{0},z_{1}-z_{0})+\mathcal{O}(h^{3})]\}d\eta+\mathcal{O}(h^{4})\\
=~ & f'(\hat{z})f'(\hat{z})\int_{-\frac{1}{2}}^{\frac{1}{2}}f(\hat{z}+\eta(z_{1}-z_{0}))d\eta+\int_{-\frac{1}{2}}^{\frac{1}{2}}[\eta^{2}f'(\hat{z})f''(\hat{z})(z_{1}-z_{0},f'(\hat{z})(z_{1}-z_{0}))\\
& +\frac{\eta^{2}}{2!}f'(\hat{z})f'''(\hat{z})(z_{1}-z_{0},z_{1}-z_{0},f(\hat{z}))+\eta^{2}f''(\hat{z})(z_{1}-z_{0},f'(\hat{z})f'(\hat{z})(z_{1}-z_{0}))\\
& +\eta^{2}f''(\hat{z})(z_{1}-z_{0},f''(\hat{z})(z_{1}-z_{0},f(\hat{z})))+\frac{\eta^{2}}{2!}f'''(\hat{z})(z_{1}-z_{0},z_{1}-z_{0},f'(\hat{z})f(\hat{z}))]d\eta+\mathcal{O}(h^{4})\\
=~ & f'f'F+\frac{h^{2}}{12}[f'f''(f'f,F)+f''(f'f'f,F)+f''(f''(f,f),F)+\frac{1}{2}f'f'''(f,f,F)\\
& +\frac{1}{2}f'''(f'f,f,F)] +\mathcal{O}(h^{4}).
\end{align*}
Considering $|\tau|=5$, we have
\begin{align*}
& \int_{-\frac{1}{2}}^{\frac{1}{2}}F_{f}(\tree58)(\hat{z}+\eta(z_{1}-z_{0}))d\eta\\
=~ & \int_{-\frac{1}{2}}^{\frac{1}{2}}f'(\hat{z}+\eta(z_{1}-z_{0}))f'(\hat{z}+\eta(z_{1}-z_{0}))f'(\hat{z}+\eta(z_{1}-z_{0}))f'(\hat{z}+\eta(z_{1}-z_{0}))f(\hat{z}+\eta(z_{1}-z_{0}))d\eta\\
=~ & \int_{-\frac{1}{2}}^{\frac{1}{2}}[f'(\hat{z})+\eta f''(\hat{z})(z_{1}-z_{0})+\mathcal{O}(h^{2})]\cdot [f'(\hat{z})+\eta f''(\hat{z})(z_{1}-z_{0})+\mathcal{O}(h^{2})]\\
& \cdot [f'(\hat{z})+\eta f''(\hat{z})(z_{1}-z_{0})+\mathcal{O}(h^{2})]\cdot [f'(\hat{z})+\eta f''(\hat{z})(z_{1}-z_{0})+\mathcal{O}(h^{2})]\\
& \cdot [f(\hat{z})+\eta f'(\hat{z})(z_{1}-z_{0})+\mathcal{O}(h^{2})]d\eta\\
=~ & \int_{-\frac{1}{2}}^{\frac{1}{2}}f'(\hat{z})f'(\hat{z})f'(\hat{z})f'(\hat{z})f(\hat{z}+\eta(z_{1}-z_{0}))d\eta+\mathcal{O}(h^{2})\\
=~ & f'f'f'f'F+\mathcal{O}(h^{2}) := A(\tree58)F+\mathcal{O}(h^{2}),
\end{align*}
where the coefficient matrix $A(\tree58)$ is defined by $f'f'f'f'F = A(\tree58)F$. In the same way, we can obtain
\begin{align*}
& \int_{-\frac{1}{2}}^{\frac{1}{2}}F_{f}(\tree59)(\hat{z}+\eta(z_{1}-z_{0}))d\eta=~ f''(f''(F,f),f)+\mathcal{O}(h^{2}) := A(\tree59)F+\mathcal{O}(h^{2}),\\
& \int_{-\frac{1}{2}}^{\frac{1}{2}}F_{f}(\tree57)(\hat{z}+\eta(z_{1}-z_{0}))d\eta=~ f'f'f''(F,f)+\mathcal{O}(h^{2}) := A(\tree57)F+\mathcal{O}(h^{2}),\\
& \int_{-\frac{1}{2}}^{\frac{1}{2}}F_{f}(\tree56)(\hat{z}+\eta(z_{1}-z_{0}))d\eta=~ f''(f'f'F,f)+\mathcal{O}(h^{2}) := A(\tree56)F+\mathcal{O}(h^{2}),\\
& \int_{-\frac{1}{2}}^{\frac{1}{2}}F_{f}(\tree54)(\hat{z}+\eta(z_{1}-z_{0}))d\eta=~ f'f'''(F,f,f)+\mathcal{O}(h^{2}) := A(\tree54)F+\mathcal{O}(h^{2}),\\
& \int_{-\frac{1}{2}}^{\frac{1}{2}}F_{f}(\tree52)(\hat{z}+\eta(z_{1}-z_{0}))d\eta=~ f'''(f'F,f,f)+\mathcal{O}(h^{2}) := A(\tree52)F+\mathcal{O}(h^{2}),\\
& \int_{-\frac{1}{2}}^{\frac{1}{2}}F_{f}(\tree53)(\hat{z}+\eta(z_{1}-z_{0}))d\eta=~ f''(f'F,f'f)+\mathcal{O}(h^{2}) := A(\tree53)F+\mathcal{O}(h^{2}),\\
& \int_{-\frac{1}{2}}^{\frac{1}{2}}F_{f}(\tree55)(\hat{z}+\eta(z_{1}-z_{0}))d\eta=~ f'f''(F,f'f)+\mathcal{O}(h^{2}) := A(\tree55)F+\mathcal{O}(h^{2}).
\end{align*}
So we have
\begin{align*}
& \int_{0}^{1}\check{g}(\xi z_{1}+(1-\xi)z_{0}))d\xi=~\int_{-\frac{1}{2}}^{\frac{1}{2}}\check{g}(\hat{z}+\eta(z_{1}-z_{0}))d\eta\\
=~ & \{I-\frac{h^{2}}{12}A(\tree32)-\frac{h^{4}}{720}[5A(\tree55)+5A(\tree56)+5A(\tree59)+\frac{5}{2}A(\tree54)+\frac{5}{2}A(\tree52)]\\
   & +\frac{h^{4}}{720}[6A(\tree58)+4A(\tree59)+A(\tree57)+4A(\tree56)+A(\tree54)+A(\tree52)+3A(\tree53)+8A(\tree55)]\}F+\mathcal{O}(h^{6})\\
=~ & \{I-\frac{h^{2}}{12}A(\tree32)+\frac{h^{4}}{720}[6A(\tree58)-A(\tree59)+A(\tree57)-A(\tree56)-\frac{3}{2}A(\tree54)-\frac{3}{2}A(\tree52)\\\\
   & +3A(\tree53)+3A(\tree55)]\}F+\mathcal{O}(h^{6}).
\end{align*}
We obtain the sixth order AVF method
\begin{align}\label{eq:4s10}
 \Phi_{h}(z_{0})=z_{1}= & c(\emptyset)z_{0}+(\sum_{\tau\in \mathcal{T}}h^{|\tau|}c(\tau)A(\tau))\int_{0}^{1}f(\xi z_{1}+(1-\xi)z_{0}))d\xi\nonumber\\
 = &z_{0}+\{hI-\frac{h^{3}}{12}A(\tree32)+\frac{h^{5}}{720}[6A(\tree58)-A(\tree59)+A(\tree57)-A(\tree56)-\frac{3}{2}A(\tree54)\nonumber\\ \nonumber\\
   & -\frac{3}{2}A(\tree52)+3A(\tree53)+3A(\tree55)]\}\int_{0}^{1}f(\xi z_{1}+(1-\xi)z_{0}))d\xi,
\end{align}
where $A(\tau)=A(\tau)(\frac{z_{1}+z_{0}}{2})$ are coefficient matrices of $F=\int_{0}^{1}f(\xi z_{1}+(1-\xi)z_{0}))d\xi$.

\begin{theorem} The sixth order AVF method $\Phi_{h}(z_{0})$ satisfies
\begin{equation*}
||\Phi_{h}(z_{0})-\varphi_{h}^{f}(z_{0})||=\mathcal{O}(h^{7}).
\end{equation*}
\end{theorem}
\begin{proof} Now we have
\begin{equation*}
\Phi_{h}^{\check{g}}(z_{0})=\Phi_{h}(z_{0})+\mathcal{O}(h^{7}), ~and\quad \Phi_{h}^{g}(z_{0})=\varphi_{h}^{f}(z_{0}),
\end{equation*}
so we can obtain
\begin{align*}
& ||\Phi_{h}(z_{0})-\varphi_{h}^{f}(z_{0})|| = ||\Phi_{h}(z_{0})-\Phi_{h}^{g}(z_{0})||= ||(\Phi_{h}^{\check{g}}(z_{0})+\mathcal{O}(h^{7}))-\Phi_{h}^{g}(z_{0})||\\
= & ||(z_{1}-z_{0})+h\int_{-1}^{1}[\check{g}(\xi z_{1}+(1-\xi)z_{0}))-g(\xi z_{1}+(1-\xi)z_{0}))]d\xi+\mathcal{O}(h^{7})||\\
= & \mathcal{O}(h^{7}).
\end{align*}
The proof is complete.
\end{proof}

\begin{theorem} If \eqref{eq:4s1} is a Hamiltonian system,
the sixth order AVF method $\Phi_{h}(z_{0})$ can preserve the discrete energy of it, i.e.
\begin{equation*}
\frac{1}{h}(H(z_{1})-H(z_{0}))=0.
\end{equation*}
\end{theorem}
\begin{proof} \eqref{eq:4s1} can be rewritten as
\begin{equation*}
\dot{z}=f(z)=S\nabla H(z),
\end{equation*}
where $S$ denotes an arbitrary constant skew-symmetric matrix and $H$ denotes the Hamiltonian. So the sixth order AVF method \eqref{eq:4s10} can be rewritten as
\begin{equation}\label{eq:4s11}
\frac{z_{1}-z_{0}}{h}=~ \tilde{S}\int_{0}^{1}\nabla H(\xi z_{1}+(1-\xi)z_{0}))d\xi,
\end{equation}
where
\begin{align*}
\tilde{S}=~ & \{I-\frac{h^{2}}{12}A(\tree32)+\frac{h^{4}}{720}[6A(\tree58)-A(\tree59)+A(\tree57)-A(\tree56)\\\\
   &-\frac{3}{2}A(\tree54) -\frac{3}{2}A(\tree52)+3A(\tree53)+3A(\tree55)]\}S,
\end{align*}
is a skew-symmetric matrix. It is given by
\begin{align*}
& IS=S,\quad A(\tree32)S=S\mathcal{H}S\mathcal{H}S,\\\\
& A(\tree58)S=S\mathcal{H}S\mathcal{H}S\mathcal{H}S\mathcal{H}S,\quad A(\tree59)S=S\mathcal{T}S\mathcal{T}S,\\\\
& A(\tree57)S=S\mathcal{H}S\mathcal{H}S\mathcal{T}S,\quad A(\tree56)S=S\mathcal{T}S\mathcal{H}S\mathcal{H}S,\\\\
& A(\tree54)S=S\mathcal{H}S\mathcal{L}S,\quad A(\tree52)S=S\mathcal{L}S\mathcal{H}S,\\\\
& A(\tree53)S=S\mathcal{R}S\mathcal{H}S,\quad A(\tree55)S=S\mathcal{H}S\mathcal{R}S,
\end{align*}
with the symmetric matrices $\mathcal{H}(z)$, $\mathcal{T}(z)$, $\mathcal{L}(z)$ and $\mathcal{R}(z)$ being given by
\begin{align*}
& \mathcal{H}_{ij}:=\frac{\partial^{2}H}{\partial z_{i}\partial z_{j}},\quad \mathcal{T}_{ij}:=\frac{\partial^{3}H}{\partial z_{i}\partial z_{j}\partial z_{k}}S^{kl}\frac{\partial H}{\partial z_{l}},\\
& \mathcal{L}_{ij}:=\frac{\partial^{4}H}{\partial z_{i}\partial z_{j}\partial z_{k}\partial z_{l}}S^{km}\frac{\partial H}{\partial z_{m}}S^{ln}\frac{\partial H}{\partial z_{n}},\quad \mathcal{R}_{ij}:=\frac{\partial^{3}H}{\partial z_{i}\partial z_{j}\partial z_{k}}S^{kl}\frac{\partial^{2} H}{\partial z_{l}\partial z_{m}}S^{mn}\frac{\partial H}{\partial z_{n}},
\end{align*}
and
\begin{align*}
& (S\mathcal{H}S\mathcal{H}S\mathcal{T}S-S\mathcal{T}S\mathcal{H}S\mathcal{H}S)^{T}=S\mathcal{T}S\mathcal{H}S\mathcal{H}S-S\mathcal{H}S\mathcal{H}S\mathcal{T}S
=-(S\mathcal{H}S\mathcal{H}S\mathcal{T}S-S\mathcal{T}S\mathcal{H}S\mathcal{H}S),\\\\
& (S\mathcal{H}S\mathcal{L}S+S\mathcal{L}S\mathcal{H}S)^{T}=-S\mathcal{L}S\mathcal{H}S-S\mathcal{H}S\mathcal{L}S
=-(S\mathcal{H}S\mathcal{L}S+S\mathcal{L}S\mathcal{H}S),\\\\
& (S\mathcal{R}S\mathcal{H}S+S\mathcal{H}S\mathcal{R}S)^{T}=-S\mathcal{H}S\mathcal{R}S-S\mathcal{R}S\mathcal{H}S
=-(S\mathcal{R}S\mathcal{H}S+S\mathcal{H}S\mathcal{R}S).
\end{align*}
We can obtain
\begin{eqnarray*}
& & \frac{1}{h}(H(z_{1})-H(z_{0}))=\frac{1}{h}\int_{0}^{1}\frac{d}{d\xi}H(\xi z_{1}+(1-\xi)z_{0})d\xi\\
& = & (\int_{0}^{1}\nabla H(\xi z_{1}+(1-\xi)z_{0})d\xi)^{T}(\frac{z_{1}-z_{0}}{h})\\
& = & (\int_{0}^{1}\nabla H(\xi z_{1}+(1-\xi)z_{0})d\xi)^{T} \tilde{S}\int_{0}^{1}\nabla H(\xi z_{1}+(1-\xi)z_{0})d\xi=0.
\end{eqnarray*}
It follows that the Hamiltonian $H$ is conserved at every time step.
\end{proof}

\paragraph{Remark 1} In the same way, we can also obtain a fifth order method
\begin{align}\label{eq:4s12}
 & \tilde{\Phi}_{h}(z_{0})=z_{1}=~ z_{0}+\{hI-\frac{h^{3}}{12}A(\tree32)-\frac{h^{4}}{24}[A(\tree42)+A(\tree43)]+\frac{h^{5}}{720}[6A(\tree58)-16A(\tree59)+A(\tree57)\\ \nonumber\\
   &-16A(\tree56)-9A(\tree54) -9A(\tree52)+3A(\tree53)-12A(\tree55)]\}\int_{0}^{1}f(\xi z_{1}+(1-\xi)z_{0}))d\xi,\nonumber
\end{align}
where $A(\tau)=A(\tau)(z_{0})$ are coefficient matrices of $F=\int_{0}^{1}f(\xi z_{1}+(1-\xi)z_{0}))d\xi$. This method is of order five, but it can not preserve the Hamiltonian, because
 $\tilde{S}$ which is the corresponding total coefficient matrix of $F$ turns out to be not skew-symmetric when expanding $F_{f}(\tau)(\xi z_{1}+(1-\xi)z_{0})$ in a Taylor series about $z=z_{0}$.

\paragraph{Remark 2} Omitting the items containing $h^{5}$ in \eqref{eq:4s10}, the method
\begin{equation*}
z_{1}=z_{0}+(hI-\frac{h^{3}}{12}A(\tree32))\int_{0}^{1}f(\xi z_{1}+(1-\xi)z_{0}))d\xi
\end{equation*}
is the fourth order method \eqref{eq:2s4}, where $A(\tree32)$ is defined in \eqref{eq:4s10}.

\begin{table}
\caption{Coefficients $\sigma(\tau)$, $\gamma(\tau)$, $a(\tau)$, $b(\tau)$ and $c(\tau)$ for trees of order $\leqslant 5$}
\label{tab:3}       
\begin{tabular}{llllllllll}
\hline\noalign{\smallskip}
&      &     &     &     &     &     &     &     &       \\
&      &     &     &     &     &     &     &     &       \\
$\tau$  &  $\emptyset$  & \tree11   &  \tree21  &  \tree31  &  \tree32  &  \tree41  &  \tree42  &  \tree43  &  \tree44  \\   \hline\noalign{\smallskip}
$\sigma(\tau)$ &    &  1  &  1  &  2  &  1  &  6  &   1 &  2  &   1   \\
$\gamma(\tau)$ &    &  1  &  2  &  3  &   6 &  4  &  8  &  12  &   24   \\
$a(\tau)$      &  1  &  1  &  1/2  &  1/3  &  1/4  &  1/4  &  1/6  &  1/6  &  1/8    \\
$b(\tau)$      &  0  &  1  &  0  &  0  &  -1/12  & 0   &  0  &  0  &  0    \\
$c(\tau)$      &  1  &  1  &  0  &  0  &  -1/12  & 0   &  0  &  0  &  0    \\   \hline\noalign{\smallskip}\hline\noalign{\smallskip}
     &      &    &    &    &    &    &    &    &      \\
     &      &    &    &    &    &    &    &    &      \\
$\tau$  &  \tree51  &  \tree52  &  \tree53  &  \tree54  &  \tree55  &  \tree56  &  \tree57  &  \tree58  &  \tree59     \\  \hline\noalign{\smallskip}
$\sigma(\tau)$ &  24  &  2  &  2  &  6  &  1  &  1  &  2  &  1  &  2    \\
$\gamma(\tau)$ &  5  &   10 &  20  & 20   &  40  &  30  & 60   &  120  &  15    \\
$a(\tau)$      &  1/5  &  1/8  & 1/12   &  1/8  &  1/12  &  1/12  &  1/12  &  1/16  &  1/9    \\
$b(\tau)$      & 0   &  1/360  & 1/120   &  1/120  &  1/90  &  1/180  &  1/360  &  1/120  &   1/90   \\
$c(\tau)$      & 0   &  -1/480 & 1/240   &  -1/480 &  1/240 &  -1/720 &  1/720  &  1/120  &   -1/720    \\ \hline\noalign{\smallskip}
\end{tabular}
\end{table}

\section{Numerical simulations}
\label{sec:5}

We here report a few numerical experiments, in order to illustrate our results presented in the previous section.

The relative energy error at $t=t_{j}$ is defined as
\begin{equation*}
RH_{j}=\frac{|H_{j}-H_{0}|}{|H_{0}|},
\end{equation*}
where $H_{j}$ denotes the Hamiltonian at $t=t_{j}$, $j=0,1,\ldots$\\
The solution error at $t=t_{j}$ is defined as
\begin{equation*}
error_{j}(h)=||z_{j}-z(t_{j})||_{\infty},
\end{equation*}
where $h$ is the time step.\\
We define
\begin{equation*}
order=log_{2}(\frac{error_{j}(h)}{error_{j}(h/2)}),
\end{equation*}
and recall that for a $p$-order accurate scheme
\begin{equation*}
\frac{error_{j}(h)}{error_{j}(h/2)}\approx 2^{p}(h\rightarrow 0).
\end{equation*}
\subsection{Numerical example 1}

First, we consider the initial value problem
\begin{align}\label{eq:5s1}
& \dot{z}=-10(z-1)^{2},\\
& z(0)=2,\nonumber
\end{align}
which has the exact solution $z=1/(1+10t)+1$. We can obtain the sixth order AVF method of \eqref{eq:5s1}
\begin{align}\label{eq:5s2}
z_{1}=~ & z_{0}+[h-\frac{h^{3}}{12}400(\frac{z_{1}+z_{0}}{2}-1)^{2}+\frac{h^{5}}{720}(14\times10^{5})(\frac{z_{1}+z_{0}}{2}-1)^{4}]\cdot\nonumber\\
& [-\frac{10}{3}((z_{1})^{2}+z_{1}z_{0}+(z_{0})^{2})+10(z_{1}+z_{0})-10].
\end{align}
The second, fourth and sixth order AVF methods and the second order trapezoidal method are applied to solve the problem \eqref{eq:5s1} in the interval $[0, 5]$ with different time steps. Table \ref{tab:4} shows the comparison of the solution errors of the four methods. We can see that the solution
errors of the second order AVF method and the second order trapezoidal method have little difference, and the three AVF methods are of order two, four and six respectively when solving the example 1.

\begin{table}
\caption{The convergence order of the second, fourth and sixth order AVF methods and the second order trapezoidal method}
\label{tab:4}       
\begin{tabular} {l  l  l  l  l} \hline\noalign{\smallskip}
\quad $t=5$   ~~~~~~~~~~~~~~~~    &$h$  ~~~~~~~~~~~~~~~~  &error    ~~~~~~~~~~~~~~~~     &order      \\   \hline\noalign{\smallskip}
second order AVF  &      0.04   &     2.0458e-5    &      -         \\
             &      0.02   &     5.0474e-6    &    2.0191   \\
             &      0.01   &     1.2578e-6    &    2.0046     \\
             &      0.005  &     3.1419e-7    &    2.0012    \\\\
fourth order AVF  &      0.04   &     6.6634e-7    &      -         \\
             &      0.02   &     4.0323e-8    &    4.0466   \\
             &      0.01   &     2.4990e-9    &    4.0122    \\
             &      0.005  &     1.5586e-10   &    4.0030   \\\\
sixth order AVF  &      0.04   &     4.2320e-8    &      -         \\
             &      0.02   &     6.6074e-10   &    6.0011  \\
             &      0.01   &     1.0325e-11   &    5.9999    \\
             &      0.005  &     1.6231e-13   &    5.9912   \\\\
trapezoidal  &      0.04   &     3.0981e-5    &      -         \\
             &      0.02   &     7.5878e-6    &    2.0296     \\
             &      0.01   &     1.8877e-6    &    2.0071     \\
             &      0.005  &     4.7135e-7    &    2.0018    \\ \hline\noalign{\smallskip}
\end{tabular}
\end{table}

\subsection{Numerical example 2}

Second, we consider a linear Hamiltonian system
\begin{equation}\label{eq:5s3}
\dot{z}=J^{-1}\nabla H, H(z)=\frac{1}{2}(p^{2}+q^{2}),
\end{equation}
\[
z=
\left(
       \begin{array}{c}
              p   \\
              q
             \end{array}
\right),
z_{0}=
\left(
       \begin{array}{c}
              1   \\
              0
             \end{array}
\right),
J=
\left(
       \begin{array}{cc}
              0   & 1  \\
              -1  & 0
             \end{array}
\right),
\]
which has the exact solution
\begin{align*}
& p(t)=cos(t),\\
& q(t)=sin(t).
\end{align*}
We can obtain the sixth order AVF method of \eqref{eq:5s3}
\[
\left(
       \begin{array}{c}
              p_{1}   \\
              q_{1}
             \end{array}
\right)=
\left(
       \begin{array}{c}
              p_{0}   \\
              q_{0}
             \end{array}
\right)+
(h+\frac{h^{3}}{12}+\frac{h^{5}}{120})
\left(
       \begin{array}{c}
              -\frac{q_{1}+q_{0}}{2}   \\
              \frac{p_{1}+p_{0}}{2}
             \end{array}
\right)
\]
We use the second, fourth and sixth
order AVF methods to resolve the
problem in the interval $[0,200]$.

\begin{figure}
  (a)\includegraphics[width=0.45\textwidth]{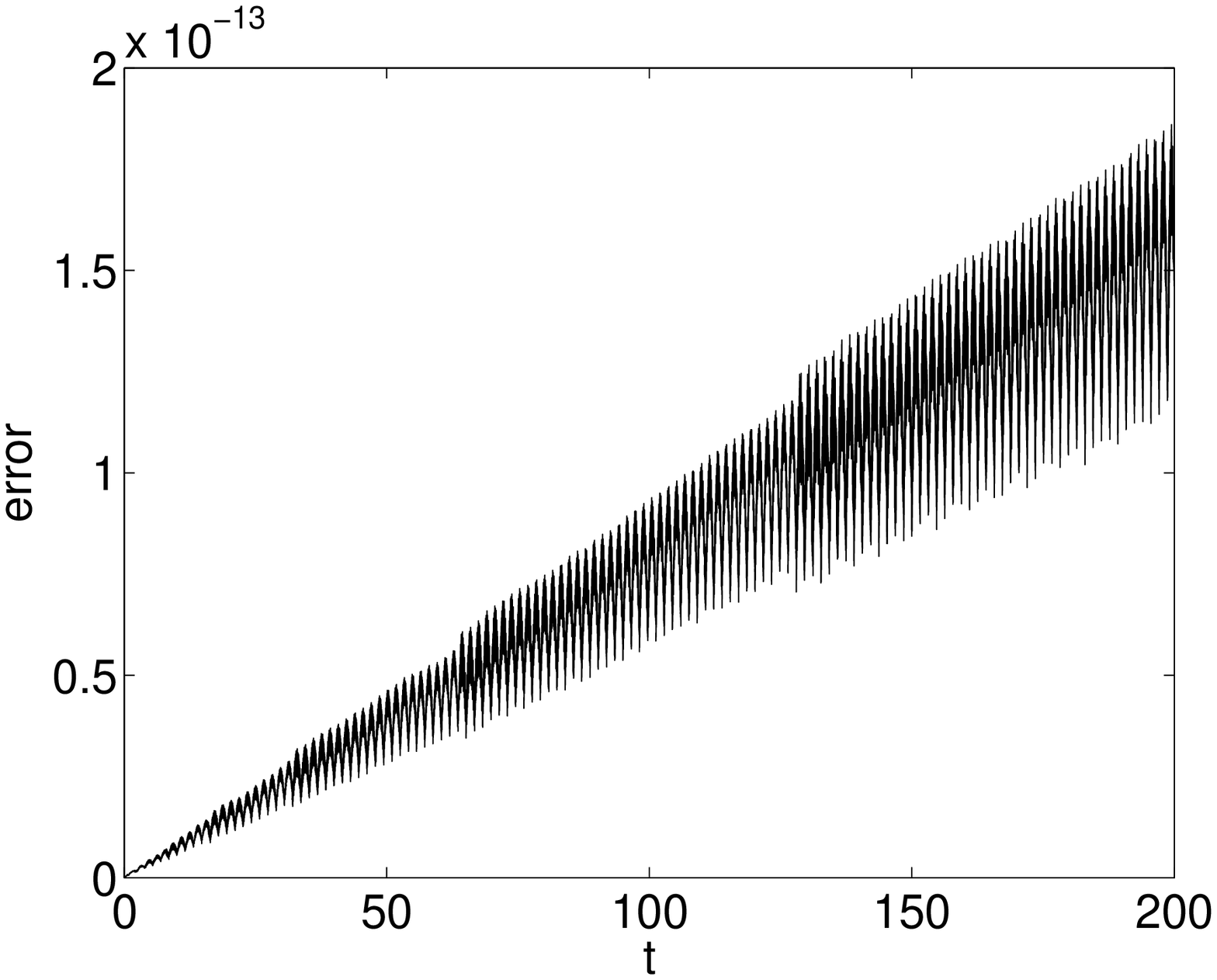}
  (b)\includegraphics[width=0.45\textwidth]{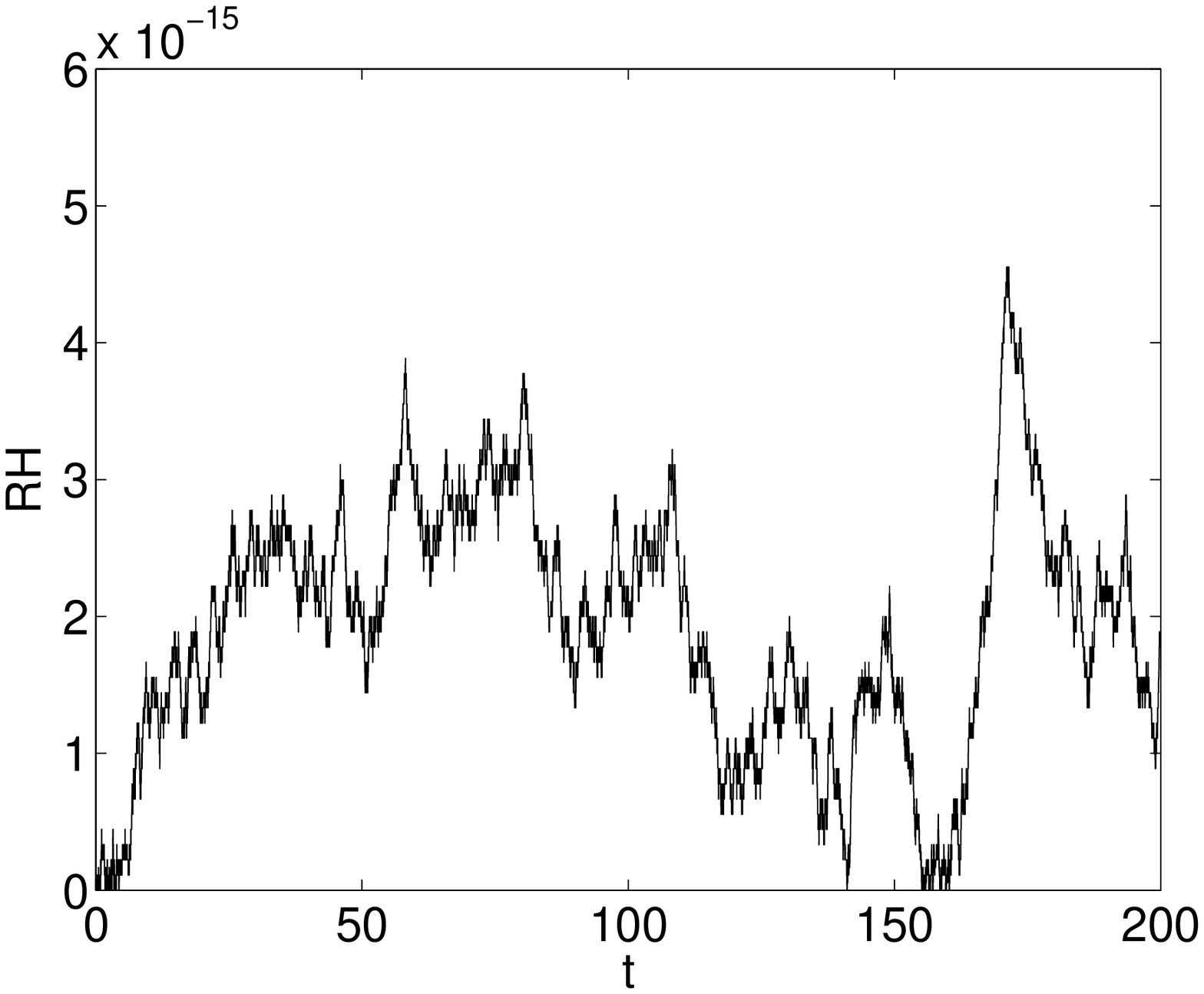}
\caption{The numerical solution of the sixth order method with $h=0.01$ from $t=0$ to $t=20000$. \textbf{a} Solution error. \textbf{b} Relative energy error}
\label{fig:1}
\end{figure}

Fig.~\ref{fig:1} shows the solution error and the relative energy error of the numerical solution of the sixth order method with time evolution in the interval $[0,200]$ respectively. From Fig.~\ref{fig:1}a we observe that the solution error of the sixth order AVF method shows a slow linear growth.
From Fig.~\ref{fig:1}b we can see that the sixth order AVF method preserves the discrete Hamiltonian very well.
So we can conclude that the sixth order AVF method has a good numerical performance for linear Hamiltonian systems in long time.

Table \ref{tab:5} shows the numerical solutions of three AVF methods with time steps $h_{1}=0.02$ and $h_{2}=0.01$. First, we can see that higher order method has higher accuracy. Second, the convergence order of the three AVF methods are 2, 4 and 6 respectively. Third, the relative energy error of the three methods is up to round-off error.

Table \ref{tab:6} shows the convergence order of the sixth order AVF methods with different time steps. We can conclude that the new method is of order six when solving the linear Hamiltonian system.

\begin{table}
\caption{The numerical solutions of three AVF methods for example 2}
\label{tab:5}       
\begin{tabular} {l l l l} \hline\noalign{\smallskip}
$t=20000$                   & second order method & fourth order method & sixth order method \\   \hline\noalign{\smallskip}
$h_{1}=0.02, error(h_{1})$  & 5.8324 $\times 10^{-3}$   & 2.3287 $\times 10^{-7}$    & 9.4389 $\times 10^{-12}$\\
$h_{2}=0.01, error(h_{2})$  & 1.4562 $\times 10^{-3}$   & 1.4554 $\times 10^{-8}$    & 1.4611 $\times 10^{-13}$\\
$order$                     & 2.0019                    & 4.0000                     & 6.0135 \\
$max(RH(h_{1}))$            & 2.5535 $\times 10^{-15}$  & 2.2204 $\times 10^{-15}$   & 3.6638 $\times 10^{-15}$\\  \hline\noalign{\smallskip}
\end{tabular}
\end{table}

\begin{table}
\caption{The convergence order of the sixth order AVF method with different time steps, $h_{1}=0.2,h_{2}=0.1,h_{3}=0.05,h_{4}=0.025,h_{5}=0.0125$}
\label{tab:6}       
\begin{tabular} {l l l l l} \hline\noalign{\smallskip}
       & $\frac{error(h_{1})}{error(h_{2})}$ & $\frac{error(h_{2})}{error(h_{3})}$ & $\frac{error(h_{3})}{error(h_{4})}$ & $\frac{error(h_{4})}{error(h_{5})}$ \\   \hline\noalign{\smallskip}
$t=1$  & 63.7145 & 63.9285 & 63.8679 & 74.4285\\
$t=2$  & 63.7132 & 63.9320 & 63.9419 & 60.7933\\
$t=3$  & 63.7142 & 63.9287 & 63.9125 & 66.0209\\
$t=4$  & 63.7138 & 63.9287 & 63.9132 & 67.6980\\
$t=5$  & 63.7141 & 63.9298 & 63.9419 & 63.3118\\
$t=200$& 63.7162 & 63.9280 & 63.9409 & 62.9838\\  \hline\noalign{\smallskip}
\end{tabular}
\end{table}

\subsection{Numerical example 3}

Third, we consider a nonlinear Hamiltonian system
\begin{equation}\label{eq:5s5}
\dot{z}=J^{-1}\nabla H, H(z)=\frac{1}{4}(p^{2}+q^{2})^{2},
\end{equation}
\[
z=
\left(
       \begin{array}{c}
              p   \\
              q
             \end{array}
\right),
z_{0}=
\left(
       \begin{array}{c}
              1   \\
              0
             \end{array}
\right),
J=
\left(
       \begin{array}{cc}
              0   & 1  \\
              -1  & 0
             \end{array}
\right),
\]
which has the exact solution
\begin{align*}
& p(t)=cos(t),\\
& q(t)=sin(t).
\end{align*}
We can obtain the sixth order AVF method of \eqref{eq:5s5}
\[
\left(
       \begin{array}{c}
              p_{1}   \\
              q_{1}
             \end{array}
\right)=
\left(
       \begin{array}{c}
              p_{0}  \\
              q_{0}
             \end{array}
\right)+\{hI-\frac{h^{3}}{12}A(\tree32)+\frac{h^{5}}{720}[6A(\tree58)-A(\tree59)+A(\tree57)-A(\tree56)-\frac{3}{2}A(\tree54)
\]
$$-\frac{3}{2}A(\tree52)+3A(\tree53)+3A(\tree55)]\}F,$$
where $A(\tau)=A(\tau)(\frac{z_{1}+z_{0}}{2})$ are coefficient matrices of $F=\int_{0}^{1}f(\xi z_{1}+(1-\xi)z_{0}))d\xi$, and
\[
F=
\left(
       \begin{array}{c}
              F^{1}   \\
              F^{2}
             \end{array}
\right)=
\left(
       \begin{array}{c}
              \int_{0}^{1}f^{1}(\xi z_{1}+(1-\xi)z_{0}))d\xi  \\
              \int_{0}^{1}f^{2}(\xi z_{1}+(1-\xi)z_{0}))d\xi
             \end{array}
\right)
\]
\[
=
\left(
       \begin{array}{c}
 -[((p_{0})^{2}+\frac{1}{3}(p_{1}-p_{0})^{2}+p_{0}(p_{1}-p_{0}))q_{0}\\
 + (\frac{1}{2}(p_{0})^{2}+\frac{1}{4}(p_{1}-p_{0})^{2}+\frac{2}{3}p_{0}(p_{1}-p_{0}))(q_{1}-q_{0})\\
 + ((q_{0})^{2}+\frac{1}{3}(q_{1}-q_{0})^{2}+q_{0}(q_{1}-q_{0}))q_{0}\\
 + (\frac{1}{2}(q_{0})^{2}+\frac{1}{4}(q_{1}-q_{0})^{2}+\frac{2}{3}q_{0}(q_{1}-q_{0}))(q_{1}-q_{0})],\\
 ((p_{0})^{2}+\frac{1}{3}(p_{1}-p_{0})^{2}+p_{0}(p_{1}-p_{0}))p_{0}\\
 + (\frac{1}{2}(p_{0})^{2}+\frac{1}{4}(p_{1}-p_{0})^{2}+\frac{2}{3}p_{0}(p_{1}-p_{0}))(p_{1}-p_{0})\\
 + ((q_{0})^{2}+\frac{1}{3}(q_{1}-q_{0})^{2}+q_{0}(q_{1}-q_{0}))p_{0}\\
 + (\frac{1}{2}(q_{0})^{2}+\frac{1}{4}(q_{1}-q_{0})^{2}+\frac{2}{3}q_{0}(q_{1}-q_{0}))(p_{1}-p_{0}).
             \end{array}
\right)
\]
We use the three AVF methods to resolve the
problem from $t=0$ to $t=200$.

\begin{figure}
  (a)\includegraphics[width=0.45\textwidth]{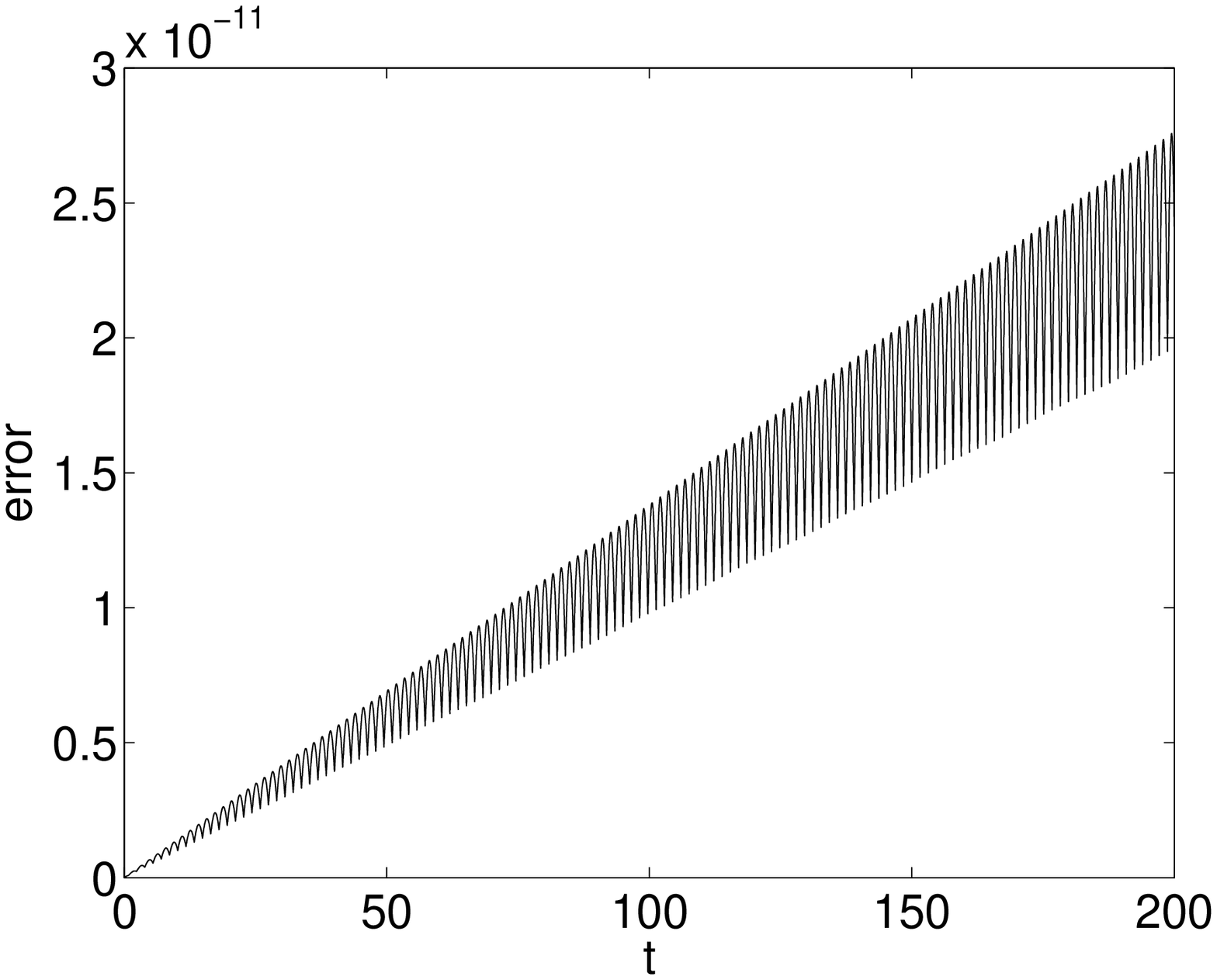}
  (b)\includegraphics[width=0.45\textwidth]{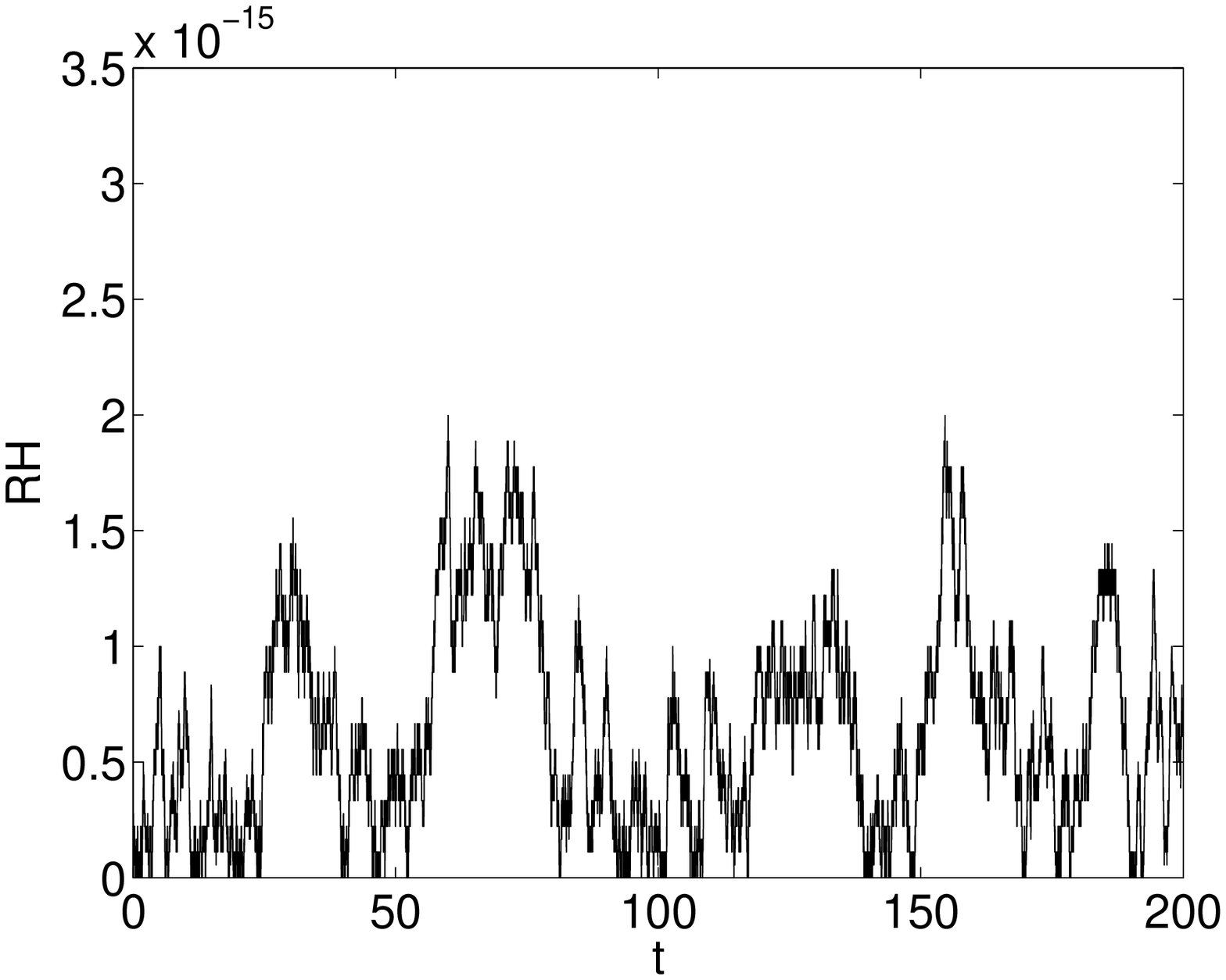}
\caption{The numerical solution of the sixth order method with $h=0.01$. \textbf{a} Solution error. \textbf{b} Relative energy error}
\label{fig:2}
\end{figure}

Fig.~\ref{fig:2} shows the solution error and the relative energy error of and the numerical solution of the sixth order method with time evolution in the interval $[0,200]$ respectively. From Fig.~\ref{fig:2}a we observe that the solution error of the sixth order AVF method shows a slow linear growth.
From Fig.~\ref{fig:2}b we can see that the sixth order AVF method preserves the discrete Hamiltonian very well.
So we can conclude that the sixth order AVF method has a good numerical performance for nonlinear Hamiltonian systems in long time.

Table \ref{tab:7} shows the numerical solutions of three AVF methods with time steps $h_{1}=0.02$ and $h_{2}=0.01$. First, we can see that higher order method has higher accuracy. Second, the convergence order of the three AVF methods are 2, 4 and 6 respectively. Third, the relative energy error of the three methods is up to round-off error.

Table \ref{tab:8} shows the convergence order of the sixth order AVF methods with different time steps. We can conclude that the new method is of order six when solving the nonlinear Hamiltonian system.

\begin{table}
\caption{The numerical solutions of three AVF methods for example 3}
\label{tab:7}       
\begin{tabular} {l l l l} \hline\noalign{\smallskip}
$t=20000$                   & second order method & fourth order method & sixth order method \\   \hline\noalign{\smallskip}
$h_{1}=0.02, error(h_{1})$  & 1.7558 $\times 10^{-2}$   & 4.5007 $\times 10^{-6}$    & 1.5495 $\times 10^{-9}$\\
$h_{2}=0.01, error(h_{2})$  & 4.3723 $\times 10^{-3}$   & 2.8137 $\times 10^{-7}$    & 2.4127 $\times 10^{-11}$\\
$order$                     & 2.0057                    & 3.9996                     & 6.0050 \\
$max(RH(h_{1}))$            & 5.2180 $\times 10^{-15}$  & 1.9984 $\times 10^{-15}$   & 1.9984 $\times 10^{-15}$\\  \hline\noalign{\smallskip}
\end{tabular}
\end{table}

\begin{table}
\caption{The convergence order of the sixth order AVF method with different time steps, $h_{1}=0.2,h_{2}=0.1,h_{3}=0.05,h_{4}=0.025,h_{5}=0.0125$}
\label{tab:8}       
\begin{tabular} {l l l l l} \hline\noalign{\smallskip}
       & $\frac{error(h_{1})}{error(h_{2})}$ & $\frac{error(h_{2})}{error(h_{3})}$ & $\frac{error(h_{3})}{error(h_{4})}$ & $\frac{error(h_{4})}{error(h_{5})}$ \\   \hline\noalign{\smallskip}
$t=1$  & 61.6174 & 63.3983 & 63.8492 & 63.9435\\
$t=2$  & 61.6178 & 63.3983 & 63.8497 & 63.9583\\
$t=3$  & 61.6174 & 63.3982 & 63.8503 & 63.9302\\
$t=4$  & 61.6167 & 63.3982 & 63.8504 & 63.9213\\
$t=5$  & 61.6179 & 63.3982 & 63.8515 & 63.9314\\
$t=200$& 61.6461 & 63.3975 & 63.9348 & 63.4203\\   \hline\noalign{\smallskip}
\end{tabular}
\end{table}

\subsection{Numerical example 4}

\begin{figure}
  (a)\includegraphics[width=0.45\textwidth]{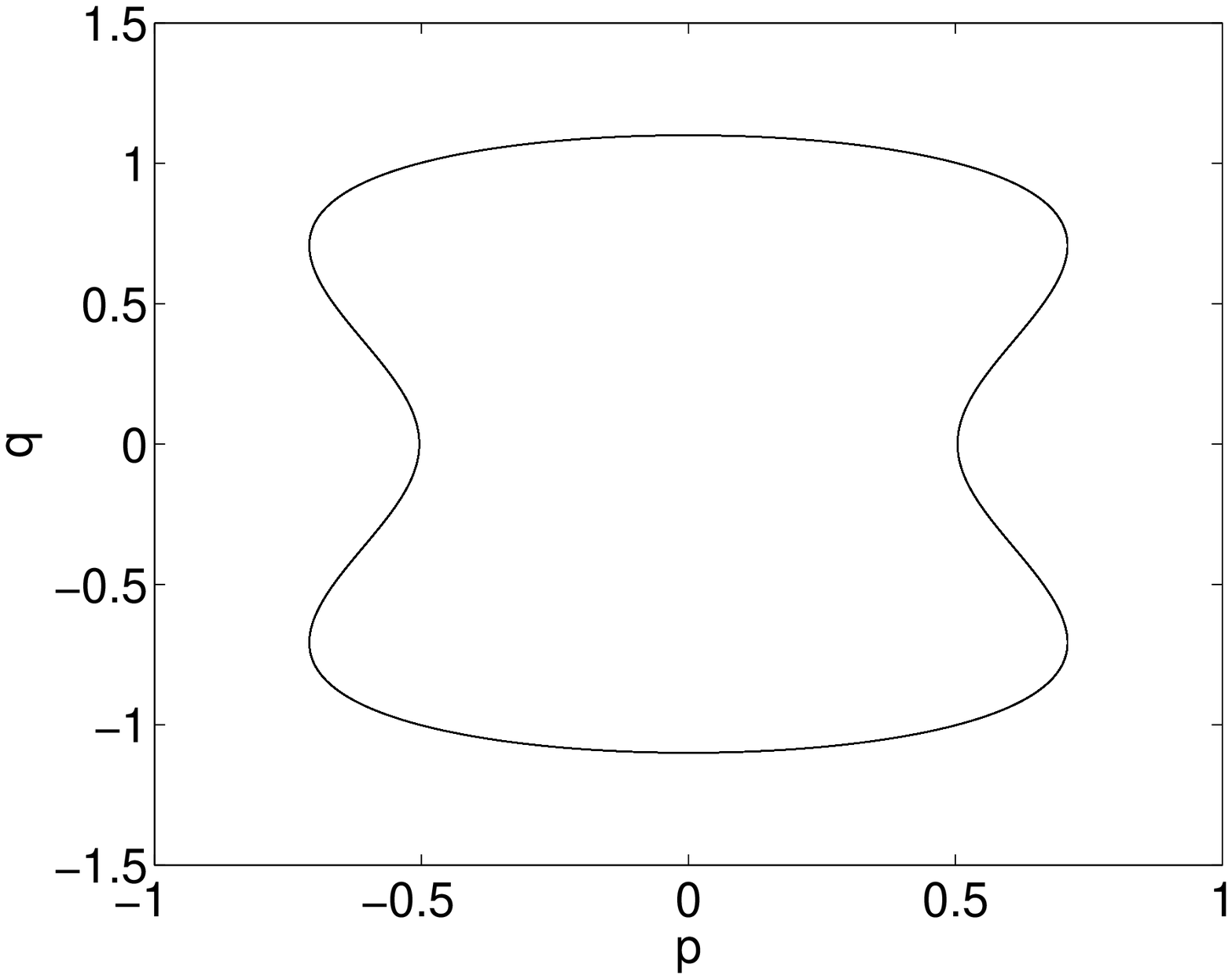}
  (b)\includegraphics[width=0.45\textwidth]{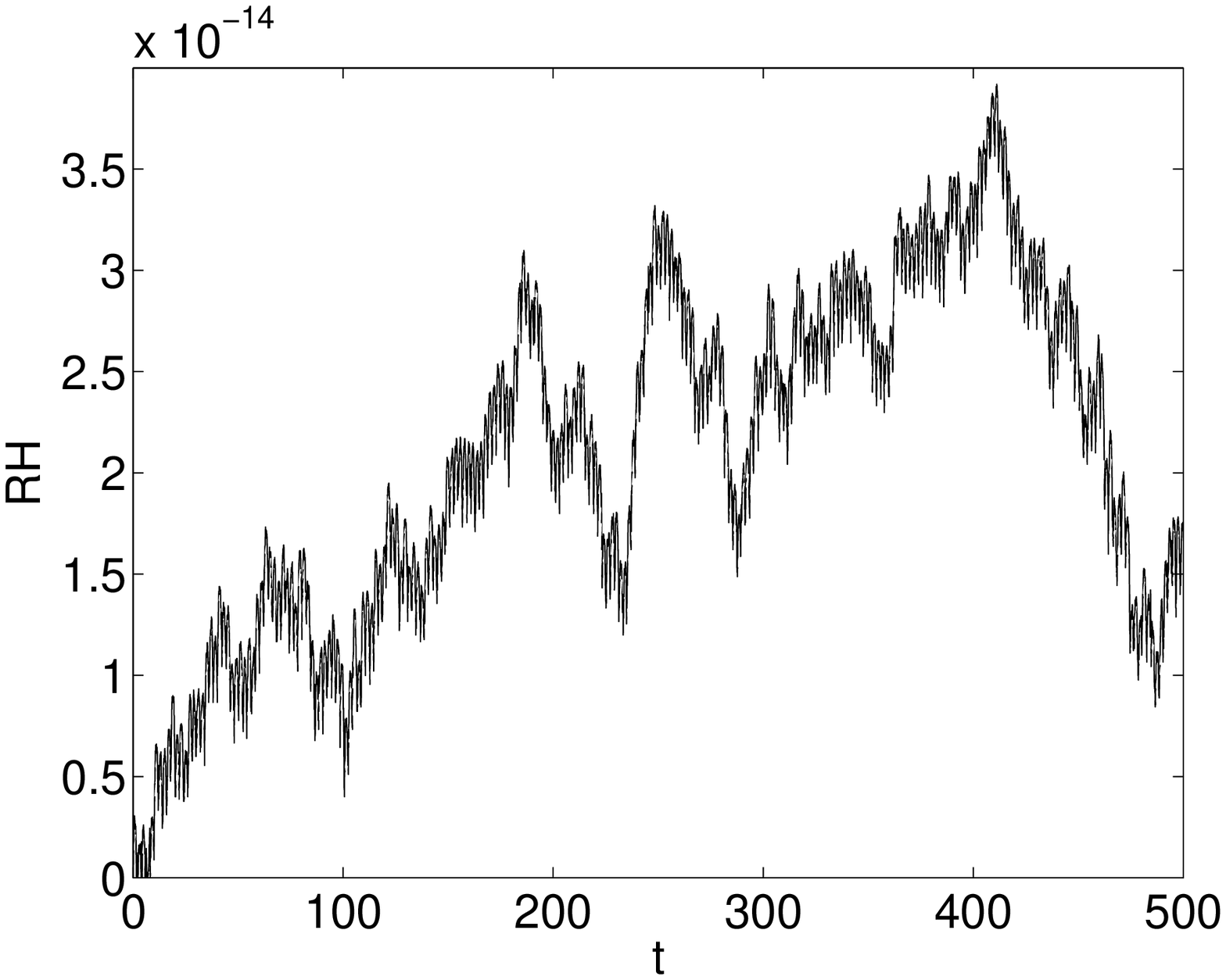}
\caption{The numerical solution of the sixth order method for the Huygens system with $h=0.005$. \textbf{a} Numerical solution. \textbf{b} Relative energy error}
\label{fig:3}
\end{figure}

Next, we consider the 2-D nonlinear Huygens system
\begin{equation}\label{eq:5s6}
\dot{z}=J^{-1}\nabla H, H(z)=p^{2}-q^{2}+q^{4},
\end{equation}
\[
z=
\left(
       \begin{array}{c}
              p   \\
              q
             \end{array}
\right),
z_{0}=
\left(
       \begin{array}{c}
              0   \\
              1.1
             \end{array}
\right),
J=
\left(
       \begin{array}{cc}
              0   & 1  \\
              -1  & 0
             \end{array}
\right).
\]
We use the sixth order AVF method to resolve the
problem from in the interval $[0,500]$. Fig.~\ref{fig:3} shows the numerical solution of the sixth order method for the Huygens system with $h=0.005$. From Fig.~\ref{fig:3}a we can see that the sixth order AVF method can solve the problem very well. From Fig.~\ref{fig:3}b we can conclude that the Hamiltonian is preserved exactly.

\subsection{Numerical example 5}

\begin{figure}
  (a)\includegraphics[width=0.45\textwidth]{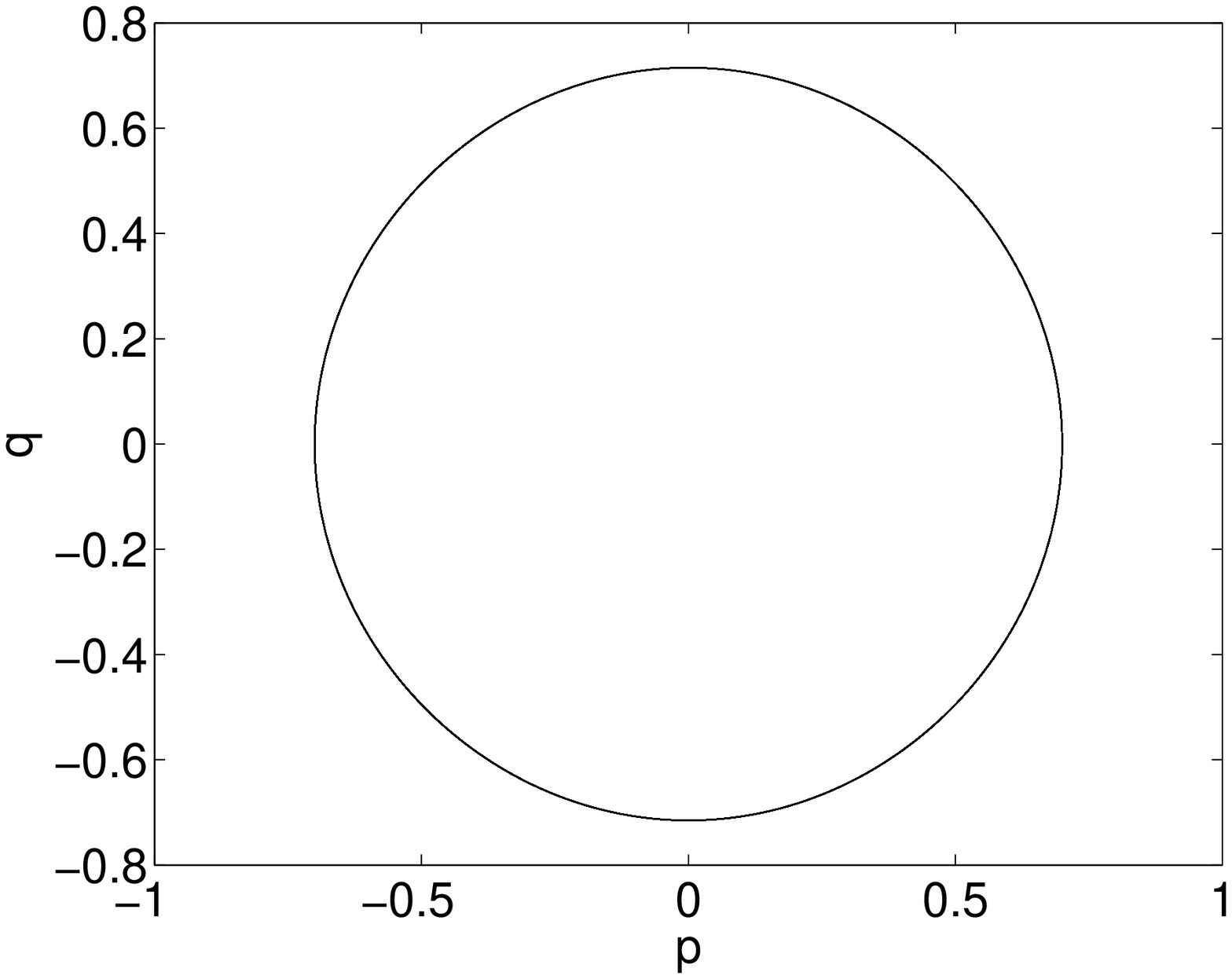}
  (b)\includegraphics[width=0.45\textwidth]{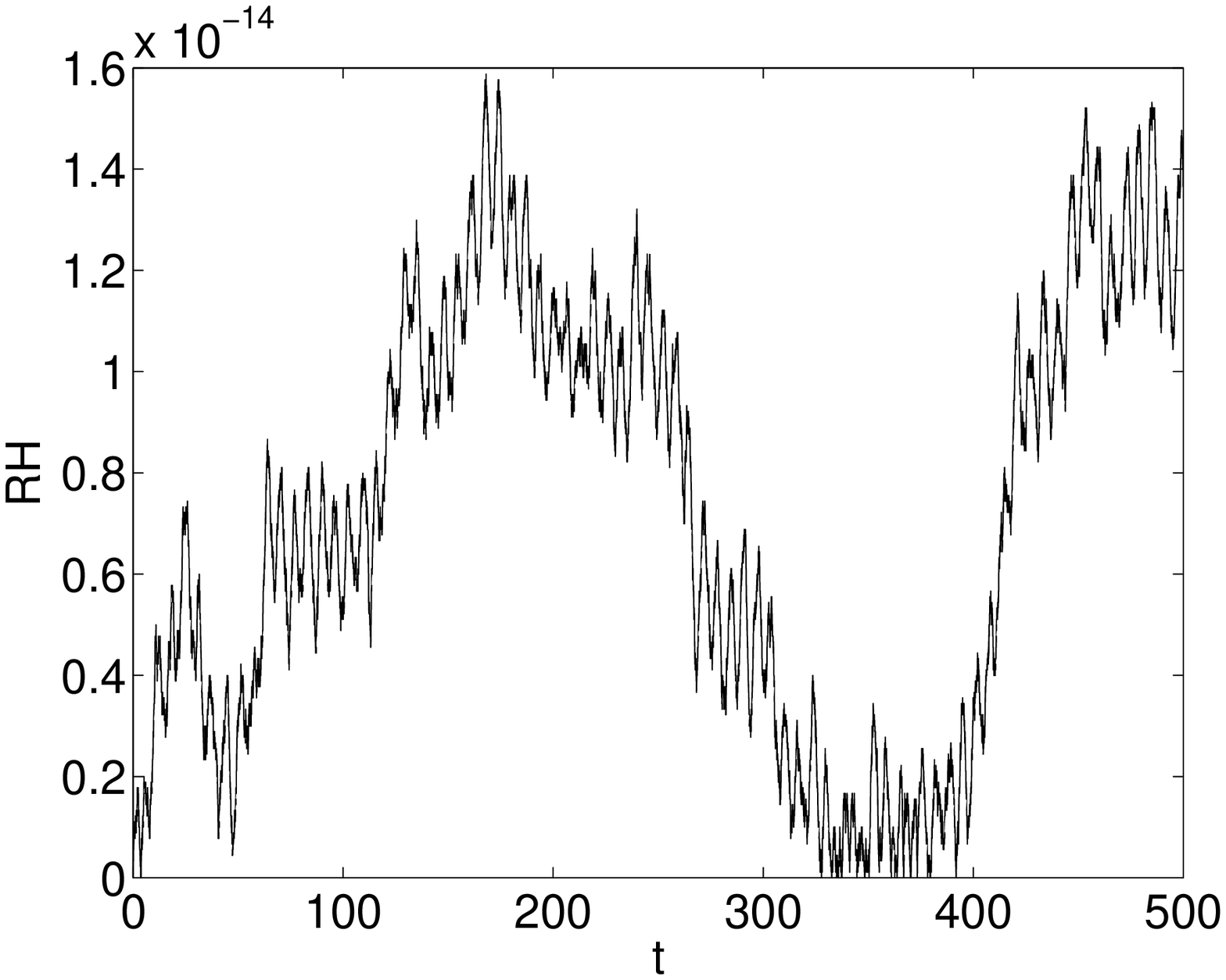}
\caption{The numerical solution of the sixth order method for the pendulum problem with $h=0.01$. \textbf{a} Numerical solution. \textbf{b} Relative energy error}
\label{fig:4}
\end{figure}

Finally, we consider the pendulum problem, which can be written as a Hamiltonian system
\begin{equation}\label{eq:5s7}
\dot{z}=J^{-1}\nabla H, H(z)=\frac{1}{2}p^{2}-cos(q),
\end{equation}
\[
z=
\left(
       \begin{array}{c}
              p   \\
              q
             \end{array}
\right),
z_{0}=
\left(
       \begin{array}{c}
              0.7   \\
              0
             \end{array}
\right),
J=
\left(
       \begin{array}{cc}
              0   & 1  \\
              -1  & 0
             \end{array}
\right).
\]
We use the sixth order AVF method to resolve the
problem in the interval $[0,500]$. Fig.~\ref{fig:4} shows the numerical solution of the sixth order method for the Huygens system with $h=0.005$. From Fig.~\ref{fig:4}a we can see that the sixth order AVF method can solve the problem very well. From Fig.~\ref{fig:4}b we can conclude that the Hamiltonian is preserved exactly.

\section{Conclusions}
\label{sec:6}

In this paper, we have proposed the concrete formulas of the substitution law for the trees of order $=5$.
Based on the new obtained substitution law, we have derived a B-series integrator extending the second order AVF method to sixth order.
This approach that we expand a low order B-series integrator into a sixth order integrator can also be used in other geometric B-series integrators naturally, easily and automatically.
We have proved that the new method is of order six and it can preserve the energy of Hamiltonian systems.
In \cite{Ref20}, Faou et al have derived the conditions a B-series method must satisfy in order to be energy-preserving.
This new method is a practical integrator of order six.
We use the sixth order AVF method to solve linear and nonlinear Hamiltonian systems to test the
accuracy and the energy-preserving ability of it. Numerical results confirm the theoretical results.

\begin{acknowledgements}
This work is supported by the Jiangsu Collaborative Innovation Center for Climate
Change, the National Natural Science Foundation of China (Grant Nos. 11271195,
41231173) and the Priority Academic Program Development of Jiangsu Higher Education Institutions.
\end{acknowledgements}

\end{document}